\documentclass[oneside,reqno,11pt]{amsart}
\usepackage[T1]{fontenc}
\usepackage[utf8]{inputenc}
\usepackage{lmodern,amsmath,amssymb,mathtools,graphicx,tikz,booktabs,array}
\usepackage[a4paper,margin=29mm]{geometry}
\usepackage[colorlinks=true,linkcolor=blue!45!black,citecolor=blue!45!black,urlcolor=blue!45!black]{hyperref}
\usepackage[nameinlink,noabbrev]{cleveref}
\usetikzlibrary{arrows.meta,calc}
\hypersetup{pdftitle={Regions of totally real arrangements of rational curves},pdfauthor={Masahiko Yoshinaga},pdfsubject={Working research draft, October 2026}}
\newtheorem{theorem}{Theorem}[section]
\newtheorem{lemma}[theorem]{Lemma}
\newtheorem{proposition}[theorem]{Proposition}
\newtheorem{corollary}[theorem]{Corollary}
\theoremstyle{definition}
\newtheorem{definition}[theorem]{Definition}
\newtheorem{example}[theorem]{Example}
\newtheorem{remark}[theorem]{Remark}

\numberwithin{equation}{section}
\newcommand{\A}{\mathcal A}
\newcommand{\Z}{\mathbb Z}
\newcommand{\Q}{\mathbb Q}
\newcommand{\R}{\mathbb R}
\newcommand{\C}{\mathbb C}
\newcommand{\F}{\mathbb F}
\newcommand{\CP}{\mathbb {CP}}
\newcommand{\RP}{\mathbb {RP}}
\newcommand{\PP}{\mathbb P}
\newcommand{\VG}{\mathcal{VG}}
\newcommand{\Fil}{\mathsf F}
\newcommand{\grVG}{\mathsf {VG}}
\newcommand{\BM}{\mathrm{BM}}
\newcommand{\gr}{\operatorname{gr}}
\newcommand{\ch}{\operatorname{ch}}
\newcommand{\bch}{\operatorname{bch}}
\newcommand{\crk}{\operatorname{cr}}
\newcommand{\can}{\operatorname{can}}
\DeclareMathOperator{\rank}{rank}
\DeclareMathOperator{\Sing}{Sing}
\DeclareMathOperator{\im}{im}
\newcommand{\UR}{U_{\R}}
\newcommand{\UC}{U_{\C}}
\newcommand{\CR}{C_{\R}}
\newcommand{\CC}{C_{\C}}
\newcommand{\one}{\mathbf 1}
\newcommand{\arxiv}[1]{\href{https://arxiv.org/abs/#1}{arXiv:#1}}
\newcommand{\doi}[1]{\href{https://doi.org/#1}{doi:#1}}
\title[Regions of rational curves]{Regions of totally real arrangements of rational curves}
\author{Masahiko Yoshinaga}
\address{Department of Mathematics, The University of Osaka, Osaka, Japan}
\email{yoshinaga@math.sci.osaka-u.ac.jp}
\date{31 October 2026}

\begin{document}
\begin{abstract}
We prove a plane-curve analogue of Zaslavsky's theorem for arrangements
of real rational curves for which every complex
singularity with more than one local branch is real and all of its branches
are real. If at least one component has odd degree, the number of real
regions equals the total rational Betti number of the complexified complement.
We also give a formula for bounded affine regions which records the number
of branches at infinity and the connected components of the real affine
curve. Oriented regions generate the middle Borel--Moore homology over
$\mathbb{Q}$, and they generate it over $\mathbb{Z}$ when the component degrees are coprime.
The relations are described by the circles in the real normalization and,
in an affine chart, by winding functions. These results imply equality of
combinatorial rank and real combinatorial rank in the framework of
Brown--Dupont and Sturmfels--Telen. Examples include singular rational
cubics, arrangements of conics, Chebyshev curves, and the plane model of
the $27$ lines on a split cubic surface.
\end{abstract}
\maketitle
\tableofcontents

\section{Introduction}
\label{sec:intro}
Zaslavsky's region-counting theorem relates the real geometry of a
hyperplane arrangement to the combinatorics of its intersection lattice.
Together with the Orlik--Solomon theorem, it implies that the number
of real chambers equals the total Betti number of the complexified
complement. The Varchenko--Gelfand
filtration refines this numerical identity degree by degree.

For arbitrary real algebraic hypersurfaces, the corresponding identity
fails. Already a smooth conic together with the line at infinity can have
nonreal intersection points, and the topology of the complexified complement
then contains information invisible in the real picture. In this paper, we study arrangements of rational plane curves whose
multibranch singularities and all their local branches are real.
Tangencies and singularities of individual branches are allowed.

Let $\A=\{C_0,\ldots,C_n\}$ be a nonempty collection of distinct,
geometrically irreducible real projective plane curves. Write
\[
 \CC=\bigcup_{i=0}^n C_i(\C),\quad
 \CR=\bigcup_{i=0}^n C_i(\R),\quad
 \UC=\CP^2\setminus\CC,\quad \UR=\RP^2\setminus\CR.
\]
The connected components of $\UR$ are called \emph{regions} or
\emph{chambers}, and their set is denoted by $\ch(\A)$.
The precise definition of a totally real arrangement of rational curves
is given in \cref{def:total}. Set $r=n+1$, $d_i=\deg C_i$, and
$\bar d=\gcd(d_0,\ldots,d_n)$. For $p\in\CC$, let $r_p$ be the number
of complex local branches, and set
\begin{equation}\label{eq:D}
 D=\sum_{p\in\Sing(\CC)}(r_p-1),\qquad h=D-r+1.
\end{equation}
Unibranch singularities, including cusps, contribute zero to $D$.

Our main results can be summarized as follows.
\begin{theorem}\label{thm:overview}
Let $\A$ be a totally real arrangement of rational plane curves in the
sense of \cref{def:total}.
\begin{enumerate}
\item If $\bar d$ is odd, every real region is an open disk and
\[
 \#\ch(\A)=D+1=\sum_{k=0}^2 b_k(\UC;\Q),\qquad
 (b_0,b_1,b_2)=(1,r-1,h).
\]
\item The classes of orientable real regions span
$H^{\BM}_2(\UC;\Q)$. If $\bar d=1$, all regions are orientable and their
classes generate $H^{\BM}_2(\UC;\Z)$; the group of relations is free of
rank $r$.
\item Suppose $C_0$ is the line at infinity, every $C_i$ for $i\geq1$
has exactly one point over infinity on its normalization, and
$\bigcup_{i=1}^n C_i(\R)\setminus C_0(\R)$ is connected. Write
$\bch(\A)$ for the set of bounded regions in this affine chart. Then
\[
 \#\bch(\A)=\chi(\UC).
\]
\item The pair $(\CP^2,\CC)$ has genus zero, and its combinatorial
and real combinatorial ranks, recalled in \cref{sec:positive}, satisfy
\[
 \crk_{\R}(\CP^2,\CC)=\crk(\CP^2,\CC)=h.
\]
\end{enumerate}
\end{theorem}
The first assertion is a plane-curve analogue of Zaslavsky's theorem (\cref{thm:regions}).
The odd-degree condition in (1) does not imply integral generation in (2): 
\cref{thm:integral} shows that the cokernel is isomorphic to $\Z/\bar d\Z$.

The last assertion uses the canonical-form map of Brown--Dupont
\cite{bd-pos} and the real combinatorial rank of Sturmfels--Telen
\cite{st-cub}. The genus-zero property of unions of rational plane curves
and the formula for their combinatorial rank are already contained in
\cite{bd-pos}. Our purpose is to relate these facts to real regions through
an explicit homological generation statement. We use positive geometry
in the sense of regions in positive arrangements when invoking
\cite{st-cub}; this should be distinguished from the original recursive
definition of Arkani-Hamed, Bai, and Lam \cite{ark}.

The proofs use normalization, the topology of a graph embedded in
$\RP^2$, and relative homology.

\section{Zaslavsky and Varchenko--Gelfand}
\label{sec:hyperplanes}
Let $\A$ now denote a finite arrangement of real affine hyperplanes in
$\R^\ell$, and put
\[
 \UR(\A)=\R^\ell\setminus\bigcup_{H\in\A}H,\qquad
 \UC(\A)=\C^\ell\setminus\bigcup_{H\in\A}H_{\C}.
\]
As above, $\ch(\A)$ denotes the set of connected components of
$\UR(\A)$; let $\bch(\A)$ be the subset of bounded chambers.
Write $\chi_{\A}(t)$ for its characteristic polynomial. The classical
formulas of Zaslavsky \cite{zas} and Orlik--Solomon \cite{orl-sol} give
\begin{equation}\label{eq:zas}
 \#\ch(\A)=(-1)^\ell\chi_{\A}(-1)
 =\sum_{k=0}^{\ell}b_k(\UC(\A)),
 \qquad
 \operatorname{Poin}(\UC(\A),t)=(-t)^\ell\chi_{\A}(-t^{-1}).
\end{equation}
If the normals span $(\R^\ell)^*$, so that the arrangement is essential,
then
\begin{equation}\label{eq:zas-bound}
 \#\bch(\A)=(-1)^\ell\chi_{\A}(1)
             =(-1)^\ell\chi(\UC(\A)).
\end{equation}
For a nonessential arrangement, the usual bounded-region count must be
replaced by relative boundedness in its essentialization. See
\cite{orl-ter} for a systematic account.

Choose a positive side $H^+$ of each hyperplane $H$, and denote the
negative side by $H^-$. The Heaviside function $x_H^+$ 
takes the value
$1$ on chambers contained in $H^+$ and $0$ on those contained in $H^-$.

The Varchenko--Gelfand algebra is
\[
 \VG(\A)=\operatorname{Map}(\ch(\A),\Z),
\]
with pointwise multiplication. Its filtration $\Fil^k\VG(\A)$ is spanned
by products of at most $k$ Heaviside functions, and
$\grVG^k(\A)=\Fil^k/\Fil^{k-1}$. The graded groups are free and
\begin{equation}\label{eq:VG-rank}
 \Fil^\ell\VG(\A)=\VG(\A),\qquad
 \rank\grVG^k(\A)=b_k(\UC(\A)).
\end{equation}
These statements originate in \cite{var-gel}; see also
\cite[Proposition 2.4]{yy-rec}. Over characteristic two, the associated
graded algebra is isomorphic to the Orlik--Solomon algebra. Over $\Z$,
the equality in \eqref{eq:VG-rank} is a statement about graded groups,
not such an identification of algebras.

Borel--Moore homology is the homology of the complex of locally finite
singular chains. Fix an orientation of $\R^\ell$. A chamber $S$ is a properly embedded,
oriented real $\ell$-manifold in $\UC(\A)$ and defines a Borel--Moore
class. Thus the region cycle map
\begin{equation}\label{eq:hyperplane-map}
 \gamma_{\A}:\VG(\A)\longrightarrow H^{\BM}_{\ell}(\UC(\A);\Z),
 \qquad \one_S\longmapsto[S],
\end{equation}
is defined. By Poincar\'e duality, its target is naturally isomorphic
to $H^{\ell}(\UC(\A);\Z)$.

The following consequence of the VG filtration and the chamber
realization of cohomology will serve as a model for our results.
See \cite{yos-ch} and \cite[Proposition 2.8 and Remark 2.10]{yy-rec}.

\begin{proposition}\label{prop:VG-top}
The map \eqref{eq:hyperplane-map} is surjective and has kernel
$\Fil^{\ell-1}\VG(\A)$.
\end{proposition}
\begin{proof}
A product of $k<\ell$ Heaviside functions comes from the subarrangement
consisting of the $k$ corresponding hyperplanes. That subarrangement has
rank at most $k$, and its complement has vanishing $\ell$th cohomology.
Restriction of Borel--Moore cycles to the open complement of the full
arrangement shows that this product maps to zero. Thus $\gamma_{\A}$
factors through $\grVG^\ell(\A)$.

We show that the induced map is an isomorphism by deletion and
restriction. If $H\in\A$, put $\A'=\A\setminus\{H\}$ and let $\A''$
be the arrangement induced on $H$. Strict compatibility of the
VG filtration with deletion--restriction
\cite[Proposition 2.8]{yy-rec} and the cohomological
deletion--restriction sequence give
\[
\begin{split}
0&\longrightarrow\grVG^\ell(\A')\longrightarrow\grVG^\ell(\A)
  \longrightarrow\grVG^{\ell-1}(\A'')\longrightarrow0,\\
0&\longrightarrow H^\ell(\UC(\A'))\longrightarrow H^\ell(\UC(\A))
  \longrightarrow H^{\ell-1}(\UC(\A''))\longrightarrow0.
\end{split}
\]
The first right-hand map is the jump of a chamber function across $H$ (see \cite[\S 2.3]{yy-rec} for details).
The second is the residue, or equivalently the boundary in the
Borel--Moore localization sequence. These maps agree on chamber classes
with compatible orientation conventions. The claim is proved by induction. 
\end{proof}
Let $F\subset\R^\ell$ be a generic affine hyperplane: it avoids all
zero-dimensional intersections of members of $\A$ and meets each
positive-dimensional intersection $X$ in dimension $\dim X-1$.
Then the classes
\[
 \{[S]:S\in\ch(\A),\ S\cap F=\varnothing\}
\]
form a $\Z$-basis of $H^{\BM}_{\ell}(\UC(\A);\Z)$
\cite[Section 2.3]{yos-ch}.

We now study the region cycle map for plane curves.

\section{Totally real arrangements and their topology}
\label{sec:setting}
We return to projective plane curves and the notation of
\cref{sec:intro}. All curves are reduced, and their listed components
are geometrically irreducible and defined over $\R$. A real curve is
\emph{geometrically irreducible} if its complexification is irreducible.
It is \emph{geometrically rational} if the normalization of its
complexification is isomorphic to $\PP^1_{\C}$.

\begin{definition}\label{def:total}
The arrangement $\A$ is a \emph{totally real arrangement of rational
curves} if:
\begin{enumerate}
\item the normalization $\widetilde C_i$ of each $C_i$ is isomorphic
over $\R$ to $\PP^1_{\R}$;
\item for every $p\in\CC$ with $r_p\geq2$, the point $p$ is real and every point of
$\nu^{-1}(p)$ is real, where
$\nu:\widetilde C=\coprod_i\widetilde C_i\to\CC$ is the normalization.
\end{enumerate}
\end{definition}
The second condition says precisely that each local complex branch at a
multibranch point is the complexification of a real branch. Tangencies
and arbitrary singularities of individual branches are allowed.
Nonreal unibranch singularities are allowed as well.

\begin{remark}\label{rem:split}
A smooth projective geometrically rational curve over $\R$ is
isomorphic to $\PP^1_{\R}$ if it has a real point $P$. Indeed,
Riemann--Roch over $\R$ gives $\dim_{\R}H^0(X,\mathcal O_X(P))=2$
for such a curve $X$. A nonconstant section is a real rational function
with a single simple pole at $P$, hence defines a degree-one morphism
$X\to\PP^1_{\R}$, which is an isomorphism.
Suppose now that the arrangement has at least two geometrically
rational components satisfying condition~(2). Each component meets
another one by B\'ezout's theorem. Such an intersection is a
multibranch point of the union, so condition~(2) supplies a real point
on the normalization of each component. The preceding fact then gives
condition~(1). 
\end{remark}

\subsection{Normalization and the branch incidence graph}
Let $P=\{p\in\CC:r_p\geq2\}$. Define a finite bipartite multigraph
$\Gamma_C$ with one vertex $v_i$ for each irreducible component, one
vertex $v_p$ for each $p\in P$, and one edge from $v_i$ to $v_p$ for
each point of $\widetilde C_i$ lying over $p$. Multiple edges must be
retained. B\'ezout's theorem makes this graph connected; for a single
component and $P=\varnothing$, it is one vertex. A self-node of $C_i$
has two distinct preimages on $\widetilde C_i$, so it gives two edges
joining $v_i$ to the same vertex $v_p$. These two edges form a cycle;
thus self-intersections are included in $\Gamma_C$.

The normalization map realizes the underlying topological space of
$\CC$ as the quotient of $r$ spheres by the finite identifications in
its fibers. 
A fiber consisting
of one point does not change the topology. Consequently
\begin{equation}\label{eq:complex-curve-homology}
 H_2(\CC;\Z)=\Z^r,\quad H_1(\CC;\Z)\simeq H_1(\Gamma_C;\Z)=\Z^h,
 \quad \chi(\CC)=2r-D,
\end{equation}
where
\[
 b_1(\Gamma_C)=\sum_{p\in P}r_p-(r+|P|)+1=D-r+1=h.
\]
The graph $\Gamma_C$ records these point identifications. Since each
normalization component has vanishing first homology,
\[
 H_1(\CC;\Z)\simeq
 \ker\bigl(\Z^{E(\Gamma_C)}\xrightarrow{\partial_\Gamma}
                 \Z^{V(\Gamma_C)}\bigr),
\]
where $\partial_\Gamma$ is the incidence map of the graph.
The finite point identifications preserve second homology and reduce
Euler characteristic by $D$.

Total reality also implies
\begin{equation}\label{eq:real-normalization}
 \nu(\widetilde C(\R))=\CR.
\end{equation}
For a real point $p$ with $r_p=1$, conjugation fixes its unique
preimage under $\nu$. If $r_p\geq2$, every point of $\nu^{-1}(p)$
is real by condition~(2) of \cref{def:total}. Therefore
$\CR$ is obtained from $r$ circles by the same finite identifications.
It is a connected embedded finite graph, up to subdivision at unibranch
singularities, and
\begin{equation}\label{eq:real-curve-homology}
 \chi(\CR)=-D,\qquad H_1(\CR;\Z)=\Z^{D+1}.
\end{equation}
Write $\nu_i:\widetilde C_i\to C_i$ for the component normalization
maps. Choose orientations of the circles $\widetilde C_i(\R)$ and put
$\alpha_i=\nu_{i*}[\widetilde C_i(\R)]$. If $j:\CR\hookrightarrow\CC$
is the inclusion, then there is a short exact sequence
\begin{equation}\label{eq:normalization-sequence}
 0\longrightarrow L:=\bigoplus_{i=0}^n\Z\alpha_i
 \longrightarrow H_1(\CR;\Z)
 \xrightarrow{j_*}H_1(\CC;\Z)\longrightarrow0.
\end{equation}
To prove this, regard each $\widetilde C_i(\R)$ as the equator in
$\widetilde C_i(\C)\simeq S^2$. All nontrivial identifications under
$\nu$ occur on these equators. Thus $\CC$ is obtained from $\CR$ by
attaching two disks for each component, with attaching map
$\nu_i|_{\widetilde C_i(\R)}$. Their cellular boundaries represent
$\alpha_i$ and $-\alpha_i$. No new one-cells are added, so $j_*$ is
surjective and its kernel is generated by the $\alpha_i$. To see that
these classes are independent, subdivide $\CR$ at its multibranch
points, adding a vertex on any component with no such point. Each
$\alpha_i$ traverses the edges of its component once, and distinct
components have disjoint edge sets. Since a graph has no two-cells,
the $\alpha_i$ are linearly independent in $H_1(\CR;\Z)$.

\subsection{Homology of the complexifed complement}
\begin{proposition}\label{prop:betti}
For an arrangement of $r$ geometrically rational projective plane
curves, without any reality assumption,
\begin{align}
 H_1(\UC;\Z)&\simeq\Z^{r-1}\oplus\Z/\bar d\Z,\label{eq:H1}\\
 H_2(\UC;\Z)&\simeq\Z^h,\label{eq:H2}\\
 H^{\BM}_2(\UC;\Z)&\simeq\Z^h\oplus\Z/\bar d\Z.\label{eq:BM2}
\end{align}
Moreover $H_k(\UC;\Z)=0$ for $k>2$ and
\begin{equation}\label{eq:chi-complement}
 \chi(\UC)=3-2r+D,\qquad \sum_k b_k(\UC;\Q)=D+1.
\end{equation}
The direct-sum decompositions are not asserted to be canonical.
\end{proposition}
\begin{proof}
Since $\CP^2$ is compact,
$H_k^{\BM}(\UC;\Z)=H_k(\CP^2,\CC;\Z)$. The relevant part of the
relative sequence is
\begin{equation}\label{eq:relative-complex}
 0\longrightarrow H_3(\CP^2,\CC)\longrightarrow\Z^r
 \xrightarrow{(d_0,\ldots,d_n)}\Z
 \longrightarrow H_2(\CP^2,\CC)
 \xrightarrow{\partial_{\C}}H_1(\CC)\longrightarrow0.
\end{equation}
In particular,
\begin{equation}\label{eq:BM-extension}
 0\longrightarrow\Z/\bar d\Z\longrightarrow H_2^{\BM}(\UC;\Z)
 \xrightarrow{\partial_{\C}}H_1(\CC;\Z)\longrightarrow0.
\end{equation}
The right-hand group is free, so this sequence splits as a sequence of
abelian groups. This proves \eqref{eq:BM2}.
Poincar\'e duality gives
$H_2^{\BM}(\UC;\Z)\simeq H^2(\UC;\Z)$ and
$H_3(\CP^2,\CC;\Z)\simeq H^1(\UC;\Z)$.

The complement of a nonempty projective plane curve is a smooth affine
surface. It has the homotopy type of a finite CW complex of dimension at
most two, by affine Lefschetz theory; see \cite{dimca}. In particular its
second homology is free. The universal coefficient theorem applied to
the groups just computed gives \eqref{eq:H1} and \eqref{eq:H2}.
Finally, additivity of Euler characteristic gives
$\chi(\UC)=3-\chi(\CC)=3-2r+D$.
\end{proof}

\subsection{The embedded real graph}
Let
\[
 \eta:H_1(\CR;\Z)\longrightarrow H_1(\RP^2;\Z)=\Z/2\Z
\]
be induced by inclusion. Identify $\Z/2\Z$ with $\F_2$ and put
\[
 \epsilon=\dim_{\F_2}\im(\eta)\in\{0,1\}.
\]
Thus $\epsilon$ records whether the image of $\eta$ is zero or all of
$\Z/2\Z$.

\begin{lemma}\label{lem:graph-complement}
For a totally real arrangement, if $\epsilon=1$, all components of
$\RP^2\setminus\CR$ are open disks. If $\epsilon=0$, exactly one
component is an open M\"obius band and all the others are open disks.
Furthermore
\begin{equation}\label{eq:parity}
 \eta(\alpha_i)=d_i\pmod2.
\end{equation}
\end{lemma}
\begin{proof}
Choose a triangulation of $\RP^2$ in which $\CR$ is a subcomplex.
A closed regular neighborhood $N$ of $\CR$ is obtained by taking small disks
about the vertices and thin bands along the edges. It is a connected
compact surface with boundary and deformation retracts onto $\CR$.
Each component of $\RP^2\setminus\CR$ is homeomorphic to the interior
of the corresponding component of $\RP^2\setminus\operatorname{int}N$.

Every boundary circle of $N$ is two-sided and separates $\RP^2$.
Each complementary component $Q$ has just one boundary circle.
To see this, suppose that $\partial Q$ contains distinct circles
$B_1$ and $B_2$. Choose $a\in B_1$ and $b\in B_2$. Since $N$ and $Q$
are connected, one can join $a$ to $b$ by an arc whose interior lies
in $\operatorname{int}N$, and return by an arc whose interior lies in
$\operatorname{int}Q$. After smoothing at the endpoints, the resulting
closed curve crosses $B_1$ transversely exactly once. This is impossible:
a closed curve crosses a separating circle an even number of times.
By the classification of compact surfaces, $Q$ is a disk or a
M\"obius band. There is at most one M\"obius band: the core curves of
two disjoint M\"obius bands would be disjoint representatives of the
nonzero class in $H_1(\RP^2;\F_2)$, whose self-intersection is $1$.

The orientation character of $N$ is the restriction of that of
$\RP^2$, and its value on a loop is given by $\eta$. If $\epsilon=1$,
$N$ contains an orientation-reversing loop. A complementary M\"obius
band would give a disjoint loop in the same nonzero mod-two homology
class, again contradicting the intersection pairing. Hence every
complementary component is a disk. If $\epsilon=0$, $N$ is orientable.
Its complementary components cannot all be disks, since capping the
boundary circles of an orientable surface with disks gives an
orientable surface. There is therefore exactly one M\"obius band.

To prove \eqref{eq:parity}, intersect $C_i$ with a generic real
projective line. All intersections are smooth and transverse; nonreal
intersection points occur in conjugate pairs. Thus the mod-two
intersection number of the real normalization cycle with a line is
$d_i$ modulo two.
\end{proof}

\section{Region-counting formulas}
\label{sec:counts}
\begin{theorem}\label{thm:regions}
For every totally real arrangement of rational curves,
\begin{equation}\label{eq:general-count}
 \#\ch(\A)=D+2-\epsilon.
\end{equation}
In particular, if some $d_i$ is odd, or equivalently if $\bar d$ is odd,
then
\[
 \#\ch(\A)=D+1=\sum_{k=0}^2b_k(\UC;\Q).
\]
More generally, this equality holds precisely when $\epsilon=1$.
\end{theorem}
\begin{proof}
We use standard relative homology and duality as in \cite{hatcher}.
All the spaces are semialgebraic and admit compatible finite
triangulations. Euler characteristic with compact supports is additive.
An open disk has compactly supported Euler characteristic $1$, while an
open M\"obius band has compactly supported Euler characteristic $0$.
Since $\chi(\RP^2)=1$ and $\chi(\CR)=-D$, \cref{lem:graph-complement}
gives $1=-D+\#\ch(\A)$ if $\epsilon=1$, and
$1=-D+\#\ch(\A)-1$ otherwise. Apply \cref{prop:betti} and
\eqref{eq:parity}.
\end{proof}

Examples with $\bar d=1$ appear in \cref{sec:examples}, with figures
in Appendix \cref{app:pictures}. 
The following examples have $\bar d=2$.

\begin{figure}[htbp]\centering
\begin{minipage}{.47\textwidth}\centering
\begin{tikzpicture}[x=1cm,y=1cm,line cap=round,line join=round]
\fill[even odd rule,blue!13] (5.5000,3.1568) -- (5.5000,0.0000) -- (2.3432,0.0000) -- (2.3129,0.1907) -- (2.2783,0.3782) -- (2.2420,0.5479) -- (2.2025,0.7067) -- (2.3210,0.6786) -- (2.4355,0.6581) -- (2.5592,0.6432) -- (2.6753,0.6359) -- (2.7998,0.6352) -- (2.9160,0.6411) -- (3.0316,0.6534) -- (3.1464,0.6721) -- (3.2600,0.6970) -- (3.3799,0.7306) -- (3.4900,0.7683) -- (3.6053,0.8153) -- (3.7104,0.8652) -- (3.8197,0.9250) -- (3.9183,0.9865) -- (4.0135,1.0534) -- (4.1364,1.1523) -- (4.2517,1.2601) -- (4.3586,1.3762) -- (4.4565,1.4999) -- (4.5449,1.6306) -- (4.6234,1.7675) -- (4.6914,1.9099) -- (4.7486,2.0569) -- (4.7947,2.2078) -- (4.8295,2.3618) -- (4.8526,2.5179) -- (4.8641,2.6753) -- (4.8638,2.8330) -- (4.8517,2.9904) -- (4.8279,3.1464) -- (4.7933,3.2975) -- (4.9505,3.2584) -- (5.1196,3.2221) -- (5.3032,3.1882) -- (5.5000,3.1568) -- cycle;
\fill[even odd rule,blue!11] (3.1568,5.5000) -- (5.5000,5.5000) -- (5.5000,3.1568) -- (5.3032,3.1882) -- (5.1196,3.2221) -- (4.9505,3.2584) -- (4.7933,3.2975) -- (4.7540,3.4273) -- (4.7044,3.5595) -- (4.6460,3.6881) -- (4.5834,3.8053) -- (4.5089,3.9252) -- (4.4316,4.0334) -- (4.3423,4.1427) -- (4.2458,4.2458) -- (4.1427,4.3423) -- (4.0334,4.4316) -- (3.9252,4.5089) -- (3.8053,4.5834) -- (3.6806,4.6497) -- (3.5595,4.7044) -- (3.4273,4.7540) -- (3.2975,4.7933) -- (3.2580,4.9521) -- (3.2217,5.1218) -- (3.1871,5.3093) -- (3.1568,5.5000) -- cycle;
\fill[even odd rule,blue!13] (0.0000,2.3432) -- (0.0000,5.5000) -- (3.1568,5.5000) -- (3.1871,5.3093) -- (3.2217,5.1218) -- (3.2580,4.9521) -- (3.2975,4.7933) -- (3.1790,4.8214) -- (3.0645,4.8419) -- (2.9408,4.8568) -- (2.8247,4.8641) -- (2.7085,4.8650) -- (2.5923,4.8595) -- (2.4766,4.8476) -- (2.3536,4.8279) -- (2.2400,4.8030) -- (2.1201,4.7694) -- (2.0100,4.7317) -- (1.8947,4.6847) -- (1.7896,4.6348) -- (1.6875,4.5792) -- (1.5817,4.5135) -- (1.4865,4.4466) -- (1.3636,4.3477) -- (1.2483,4.2399) -- (1.1414,4.1238) -- (1.0435,4.0001) -- (0.9551,3.8694) -- (0.8766,3.7325) -- (0.8086,3.5901) -- (0.7514,3.4431) -- (0.7053,3.2922) -- (0.6705,3.1382) -- (0.6474,2.9821) -- (0.6359,2.8247) -- (0.6362,2.6670) -- (0.6483,2.5096) -- (0.6721,2.3536) -- (0.7067,2.2025) -- (0.5495,2.2416) -- (0.3804,2.2779) -- (0.1968,2.3118) -- (0.0000,2.3432) -- cycle;
\fill[even odd rule,blue!11] (2.3432,0.0000) -- (0.0000,0.0000) -- (0.0000,2.3432) -- (0.1968,2.3118) -- (0.3804,2.2779) -- (0.5495,2.2416) -- (0.7067,2.2025) -- (0.7460,2.0727) -- (0.7956,1.9405) -- (0.8540,1.8119) -- (0.9166,1.6947) -- (0.9911,1.5748) -- (1.0734,1.4600) -- (1.1632,1.3511) -- (1.2542,1.2542) -- (1.3573,1.1577) -- (1.4666,1.0684) -- (1.5817,0.9865) -- (1.6947,0.9166) -- (1.8194,0.8503) -- (1.9405,0.7956) -- (2.0727,0.7460) -- (2.2025,0.7067) -- (2.2420,0.5479) -- (2.2783,0.3782) -- (2.3129,0.1907) -- (2.3432,0.0000) -- cycle;
\fill[even odd rule,blue!9] (0.7067,2.2025) -- (0.8676,2.1557) -- (1.0167,2.1046) -- (1.1530,2.0495) -- (1.2803,1.9888) -- (1.3966,1.9234) -- (1.5039,1.8522) -- (1.6021,1.7754) -- (1.6930,1.6916) -- (1.7766,1.6007) -- (1.8529,1.5029) -- (1.9238,1.3959) -- (1.9893,1.2794) -- (2.0493,1.1535) -- (2.1038,1.0187) -- (2.1547,0.8707) -- (2.2025,0.7067) -- (2.0727,0.7460) -- (1.9405,0.7956) -- (1.8194,0.8503) -- (1.6947,0.9166) -- (1.5817,0.9865) -- (1.4666,1.0684) -- (1.3573,1.1577) -- (1.2542,1.2542) -- (1.1632,1.3511) -- (1.0734,1.4600) -- (0.9911,1.5748) -- (0.9166,1.6947) -- (0.8540,1.8119) -- (0.7956,1.9405) -- (0.7460,2.0727) -- (0.7067,2.2025) -- cycle;
\fill[even odd rule,blue!5] (2.2025,0.7067) -- (2.1547,0.8707) -- (2.1038,1.0187) -- (2.0493,1.1535) -- (1.9893,1.2794) -- (1.9238,1.3959) -- (1.8529,1.5029) -- (1.7766,1.6007) -- (1.6930,1.6916) -- (1.6021,1.7754) -- (1.5039,1.8522) -- (1.3966,1.9234) -- (1.2803,1.9888) -- (1.1530,2.0495) -- (1.0167,2.1046) -- (0.8676,2.1557) -- (0.7067,2.2025) -- (0.6753,2.3373) -- (0.6524,2.4766) -- (0.6393,2.6089) -- (0.6346,2.7500) -- (0.6393,2.8911) -- (0.6534,3.0316) -- (0.6769,3.1708) -- (0.7074,3.3002) -- (0.7356,3.3958) -- (0.7712,3.4977) -- (0.8086,3.5901) -- (0.8540,3.6881) -- (0.9043,3.7836) -- (0.9551,3.8694) -- (1.0098,3.9527) -- (1.0734,4.0400) -- (1.1414,4.1238) -- (1.2080,4.1981) -- (1.2839,4.2749) -- (1.3573,4.3423) -- (1.4404,4.4112) -- (1.5201,4.4711) -- (1.6095,4.5316) -- (1.6947,4.5834) -- (1.7896,4.6348) -- (1.8795,4.6780) -- (1.9790,4.7199) -- (2.0727,4.7540) -- (2.1758,4.7860) -- (2.2723,4.8108) -- (2.3781,4.8324) -- (2.4766,4.8476) -- (2.5840,4.8589) -- (2.6836,4.8643) -- (2.7832,4.8651) -- (2.8911,4.8607) -- (2.9904,4.8517) -- (3.0973,4.8367) -- (3.1952,4.8180) -- (3.2975,4.7933) -- (3.3453,4.6293) -- (3.3962,4.4813) -- (3.4507,4.3465) -- (3.5107,4.2206) -- (3.5762,4.1041) -- (3.6471,3.9971) -- (3.7234,3.8993) -- (3.8070,3.8084) -- (3.8979,3.7246) -- (3.9961,3.6478) -- (4.1034,3.5766) -- (4.2197,3.5112) -- (4.3470,3.4505) -- (4.4833,3.3954) -- (4.6324,3.3443) -- (4.7933,3.2975) -- (4.8247,3.1627) -- (4.8476,3.0234) -- (4.8607,2.8911) -- (4.8654,2.7500) -- (4.8607,2.6089) -- (4.8466,2.4684) -- (4.8231,2.3292) -- (4.7926,2.1998) -- (4.7619,2.0963) -- (4.7288,2.0023) -- (4.6881,1.9023) -- (4.6460,1.8119) -- (4.5997,1.7236) -- (4.5449,1.6306) -- (4.4855,1.5404) -- (4.4266,1.4600) -- (4.3639,1.3825) -- (4.2920,1.3019) -- (4.2161,1.2251) -- (4.1427,1.1577) -- (4.0661,1.0939) -- (3.9799,1.0289) -- (3.8975,0.9729) -- (3.8053,0.9166) -- (3.7104,0.8652) -- (3.6205,0.8220) -- (3.5210,0.7801) -- (3.4273,0.7460) -- (3.3322,0.7163) -- (3.2277,0.6892) -- (3.1219,0.6676) -- (3.0234,0.6524) -- (2.9160,0.6411) -- (2.8164,0.6357) -- (2.7168,0.6349) -- (2.6089,0.6393) -- (2.5096,0.6483) -- (2.4027,0.6633) -- (2.3048,0.6820) -- (2.2025,0.7067) -- cycle;
\fill[even odd rule,blue!7] (4.7933,3.2975) -- (4.6324,3.3443) -- (4.4833,3.3954) -- (4.3470,3.4505) -- (4.2197,3.5112) -- (4.1034,3.5766) -- (3.9961,3.6478) -- (3.8979,3.7246) -- (3.8070,3.8084) -- (3.7234,3.8993) -- (3.6471,3.9971) -- (3.5762,4.1041) -- (3.5107,4.2206) -- (3.4507,4.3465) -- (3.3962,4.4813) -- (3.3453,4.6293) -- (3.2975,4.7933) -- (3.4273,4.7540) -- (3.5595,4.7044) -- (3.6806,4.6497) -- (3.8053,4.5834) -- (3.9252,4.5089) -- (4.0334,4.4316) -- (4.1427,4.3423) -- (4.2458,4.2458) -- (4.3423,4.1427) -- (4.4316,4.0334) -- (4.5089,3.9252) -- (4.5834,3.8053) -- (4.6460,3.6881) -- (4.7044,3.5595) -- (4.7540,3.4273) -- (4.7933,3.2975) -- cycle;
\draw[black,line width=.65pt] (0.0000,2.3432) -- (0.1605,2.3180) -- (0.3114,2.2913) -- (0.4532,2.2629) -- (0.5859,2.2331) -- (0.7113,2.2013) -- (0.8295,2.1675) -- (0.9404,2.1318) -- (1.0458,2.0936) -- (1.1440,2.0534) -- (1.2367,2.0108) -- (1.3239,1.9655) -- (1.4076,1.9167) -- (1.4857,1.8651) -- (1.5584,1.8111) -- (1.6275,1.7534) -- (1.6930,1.6916) -- (1.7548,1.6259) -- (1.8130,1.5561) -- (1.8675,1.4824) -- (1.9184,1.4048) -- (1.9657,1.3237) -- (2.0111,1.2360) -- (2.0547,1.1410) -- (2.0947,1.0428) -- (2.1329,0.9371) -- (2.1674,0.8297) -- (2.2020,0.7086) -- (2.2329,0.5866) -- (2.2620,0.4577) -- (2.2911,0.3124) -- (2.3183,0.1584) -- (2.3432,0.0000);
\draw[black,line width=.65pt] (3.1568,5.5000) -- (3.1817,5.3416) -- (3.2089,5.1876) -- (3.2380,5.0423) -- (3.2671,4.9134) -- (3.2980,4.7914) -- (3.3326,4.6703) -- (3.3671,4.5629) -- (3.4053,4.4572) -- (3.4453,4.3590) -- (3.4889,4.2640) -- (3.5343,4.1763) -- (3.5816,4.0952) -- (3.6325,4.0176) -- (3.6870,3.9439) -- (3.7452,3.8741) -- (3.8070,3.8084) -- (3.8725,3.7466) -- (3.9416,3.6889) -- (4.0143,3.6349) -- (4.0924,3.5833) -- (4.1761,3.5345) -- (4.2633,3.4892) -- (4.3560,3.4466) -- (4.4542,3.4064) -- (4.5596,3.3682) -- (4.6705,3.3325) -- (4.7887,3.2987) -- (4.9141,3.2669) -- (5.0468,3.2371) -- (5.1886,3.2087) -- (5.3395,3.1820) -- (5.5000,3.1568);
\draw[black,line width=.65pt] (4.8654,2.7500) -- (4.8630,2.8496) -- (4.8552,2.9573) -- (4.8419,3.0645) -- (4.8247,3.1627) -- (4.8030,3.2600) -- (4.7743,3.3641) -- (4.7403,3.4666) -- (4.7044,3.5595) -- (4.6605,3.6582) -- (4.6156,3.7472) -- (4.5623,3.8411) -- (4.5089,3.9252) -- (4.4516,4.0068) -- (4.3852,4.0920) -- (4.3202,4.1675) -- (4.2458,4.2458) -- (4.1675,4.3202) -- (4.0920,4.3852) -- (4.0068,4.4516) -- (3.9252,4.5089) -- (3.8411,4.5623) -- (3.7472,4.6156) -- (3.6507,4.6641) -- (3.5595,4.7044) -- (3.4666,4.7403) -- (3.3641,4.7743) -- (3.2600,4.8030) -- (3.1627,4.8247) -- (3.0645,4.8419) -- (2.9573,4.8552) -- (2.8496,4.8630) -- (2.7500,4.8654) -- (2.6504,4.8630) -- (2.5427,4.8552) -- (2.4355,4.8419) -- (2.3373,4.8247) -- (2.2400,4.8030) -- (2.1359,4.7743) -- (2.0334,4.7403) -- (1.9405,4.7044) -- (1.8418,4.6605) -- (1.7528,4.6156) -- (1.6589,4.5623) -- (1.5748,4.5089) -- (1.4865,4.4466) -- (1.4080,4.3852) -- (1.3325,4.3202) -- (1.2542,4.2458) -- (1.1854,4.1737) -- (1.1148,4.0920) -- (1.0484,4.0068) -- (0.9911,3.9252) -- (0.9377,3.8411) -- (0.8844,3.7472) -- (0.8395,3.6582) -- (0.7956,3.5595) -- (0.7569,3.4587) -- (0.7257,3.3641) -- (0.6990,3.2680) -- (0.6753,3.1627) -- (0.6581,3.0645) -- (0.6448,2.9573) -- (0.6370,2.8496) -- (0.6346,2.7500) -- (0.6370,2.6504) -- (0.6448,2.5427) -- (0.6569,2.4437) -- (0.6753,2.3373) -- (0.6990,2.2320) -- (0.7257,2.1359) -- (0.7597,2.0334) -- (0.7956,1.9405) -- (0.8395,1.8418) -- (0.8844,1.7528) -- (0.9334,1.6660) -- (0.9911,1.5748) -- (1.0484,1.4932) -- (1.1148,1.4080) -- (1.1798,1.3325) -- (1.2542,1.2542) -- (1.3325,1.1798) -- (1.4080,1.1148) -- (1.4865,1.0534) -- (1.5748,0.9911) -- (1.6589,0.9377) -- (1.7528,0.8844) -- (1.8418,0.8395) -- (1.9405,0.7956) -- (2.0413,0.7569) -- (2.1359,0.7257) -- (2.2400,0.6970) -- (2.3373,0.6753) -- (2.4437,0.6569) -- (2.5427,0.6448) -- (2.6504,0.6370) -- (2.7500,0.6346) -- (2.8579,0.6374) -- (2.9573,0.6448) -- (3.0563,0.6569) -- (3.1627,0.6753) -- (3.2600,0.6970) -- (3.3641,0.7257) -- (3.4666,0.7597) -- (3.5595,0.7956) -- (3.6582,0.8395) -- (3.7472,0.8844) -- (3.8411,0.9377) -- (3.9252,0.9911) -- (4.0068,1.0484) -- (4.0920,1.1148) -- (4.1737,1.1854) -- (4.2458,1.2542) -- (4.3146,1.3263) -- (4.3852,1.4080) -- (4.4466,1.4865) -- (4.5089,1.5748) -- (4.5623,1.6589) -- (4.6156,1.7528) -- (4.6605,1.8418) -- (4.7044,1.9405) -- (4.7431,2.0413) -- (4.7743,2.1359) -- (4.8010,2.2320) -- (4.8247,2.3373) -- (4.8431,2.4437) -- (4.8552,2.5427) -- (4.8630,2.6504) -- (4.8654,2.7500);
\draw[gray!50,thin] (0,0) rectangle (5.5,5.5);
\node[font=\small,fill=white,fill opacity=.85,text opacity=1,inner sep=1pt] at (2.7500,2.7500) {$S_{1}$};
\node[font=\small,fill=white,fill opacity=.85,text opacity=1,inner sep=1pt] at (4.0192,4.0192) {$S_{2}$};
\node[font=\small,fill=white,fill opacity=.85,text opacity=1,inner sep=1pt] at (1.4808,1.4808) {$S_{3}$};
\node[font=\small,fill=white,fill opacity=.85,text opacity=1,inner sep=1pt] at (5.0029,5.0029) {$S_{4}$};
\node[font=\small,fill=white,fill opacity=.85,text opacity=1,inner sep=1pt] at (0.4971,0.4971) {$S_{4}$};
\node[font=\small,fill=white,fill opacity=.85,text opacity=1,inner sep=1pt] at (0.5500,4.9500) {$S_{5}$};
\node[font=\small,fill=white,fill opacity=.85,text opacity=1,inner sep=1pt] at (4.9500,0.5500) {$S_{5}$};
\end{tikzpicture}
\par\smallskip (2)\end{minipage}
\hfill
\begin{minipage}{.47\textwidth}\centering
\begin{tikzpicture}[x=1cm,y=1cm,line cap=round,line join=round]
\fill[even odd rule,blue!15] (5.5000,0.0000) -- (0.0000,0.0000) -- (0.0000,5.5000) -- (5.5000,5.5000) -- (5.5000,0.0000) -- cycle (3.7749,3.7749) -- (3.7293,3.9397) -- (3.6770,4.0970) -- (3.6213,4.2383) -- (3.5570,4.3768) -- (3.4907,4.4984) -- (3.4162,4.6145) -- (3.3410,4.7133) -- (3.2581,4.8040) -- (3.1718,4.8807) -- (3.0826,4.9430) -- (2.9956,4.9884) -- (2.9025,5.0213) -- (2.8130,5.0382) -- (2.7185,5.0408) -- (2.6243,5.0278) -- (2.5353,5.0011) -- (2.4088,4.9377) -- (2.2867,4.8460) -- (2.1706,4.7271) -- (2.0584,4.5772) -- (1.9591,4.4082) -- (1.8700,4.2177) -- (1.7923,4.0082) -- (1.7251,3.7749) -- (1.5449,3.7246) -- (1.3740,3.6663) -- (1.2210,3.6035) -- (1.0733,3.5311) -- (0.9448,3.4558) -- (0.8248,3.3715) -- (0.7205,3.2822) -- (0.6328,3.1885) -- (0.5733,3.1084) -- (0.5281,3.0306) -- (0.4940,2.9514) -- (0.4703,2.8668) -- (0.4595,2.7860) -- (0.4605,2.7005) -- (0.4742,2.6153) -- (0.4989,2.5353) -- (0.5623,2.4088) -- (0.6577,2.2826) -- (0.7775,2.1667) -- (0.9228,2.0584) -- (1.0918,1.9591) -- (1.2823,1.8700) -- (1.4918,1.7923) -- (1.7251,1.7251) -- (1.7754,1.5449) -- (1.8337,1.3740) -- (1.8995,1.2143) -- (1.9689,1.0733) -- (2.0477,0.9392) -- (2.1285,0.8248) -- (2.2178,0.7205) -- (2.3115,0.6328) -- (2.3916,0.5733) -- (2.4694,0.5281) -- (2.5530,0.4925) -- (2.6332,0.4703) -- (2.7140,0.4595) -- (2.7995,0.4605) -- (2.8847,0.4742) -- (2.9647,0.4989) -- (3.0912,0.5623) -- (3.2133,0.6540) -- (3.3294,0.7729) -- (3.4416,0.9228) -- (3.5409,1.0918) -- (3.6300,1.2823) -- (3.7077,1.4918) -- (3.7749,1.7251) -- (4.0232,1.7973) -- (4.2452,1.8816) -- (4.4450,1.9788) -- (4.6249,2.0911) -- (4.7711,2.2099) -- (4.8873,2.3366) -- (4.9350,2.4045) -- (4.9719,2.4694) -- (5.0027,2.5397) -- (5.0247,2.6109) -- (5.0391,2.6960) -- (5.0405,2.7860) -- (5.0278,2.8757) -- (5.0027,2.9603) -- (4.9628,3.0480) -- (4.9122,3.1297) -- (4.8460,3.2133) -- (4.7711,3.2901) -- (4.6849,3.3640) -- (4.5826,3.4380) -- (4.4690,3.5078) -- (4.3513,3.5697) -- (4.2177,3.6300) -- (4.0824,3.6823) -- (3.9320,3.7317) -- (3.7749,3.7749) -- cycle;
\fill[even odd rule,blue!13] (3.8958,2.7500) -- (3.8884,3.0104) -- (3.8652,3.2762) -- (3.8266,3.5347) -- (3.7749,3.7749) -- (4.0530,3.6926) -- (4.2990,3.5945) -- (4.4082,3.5409) -- (4.5158,3.4804) -- (4.6145,3.4162) -- (4.6992,3.3525) -- (4.7753,3.2862) -- (4.8460,3.2133) -- (4.9062,3.1381) -- (4.9556,3.0610) -- (4.9922,2.9868) -- (5.0201,2.9070) -- (5.0366,2.8264) -- (5.0417,2.7500) -- (5.0359,2.6691) -- (5.0201,2.5930) -- (4.9922,2.5132) -- (4.9556,2.4390) -- (4.9092,2.3661) -- (4.8496,2.2908) -- (4.7795,2.2178) -- (4.6992,2.1475) -- (4.6145,2.0838) -- (4.5158,2.0196) -- (4.4082,1.9591) -- (4.2990,1.9055) -- (4.0530,1.8074) -- (3.7749,1.7251) -- (3.8266,1.9653) -- (3.8652,2.2238) -- (3.8884,2.4896) -- (3.8958,2.7500) -- cycle;
\fill[even odd rule,blue!5] (3.7749,3.7749) -- (3.8266,3.5347) -- (3.8652,3.2762) -- (3.8884,3.0104) -- (3.8958,2.7500) -- (3.8884,2.4896) -- (3.8652,2.2238) -- (3.8266,1.9653) -- (3.7749,1.7251) -- (3.5347,1.6734) -- (3.2762,1.6348) -- (3.0194,1.6121) -- (2.7500,1.6042) -- (2.4896,1.6116) -- (2.2238,1.6348) -- (1.9653,1.6734) -- (1.7251,1.7251) -- (1.6734,1.9653) -- (1.6348,2.2238) -- (1.6121,2.4806) -- (1.6042,2.7500) -- (1.6116,3.0104) -- (1.6348,3.2762) -- (1.6734,3.5347) -- (1.7251,3.7749) -- (1.9653,3.8266) -- (2.2238,3.8652) -- (2.4806,3.8879) -- (2.7500,3.8958) -- (3.0104,3.8884) -- (3.2762,3.8652) -- (3.5347,3.8266) -- (3.7749,3.7749) -- cycle;
\fill[even odd rule,blue!7] (1.7251,3.7749) -- (1.8074,4.0530) -- (1.9055,4.2990) -- (1.9591,4.4082) -- (2.0196,4.5158) -- (2.0838,4.6145) -- (2.1475,4.6992) -- (2.2178,4.7795) -- (2.2867,4.8460) -- (2.3619,4.9062) -- (2.4390,4.9556) -- (2.5176,4.9941) -- (2.5930,5.0201) -- (2.6736,5.0366) -- (2.7500,5.0417) -- (2.8309,5.0359) -- (2.9070,5.0201) -- (2.9868,4.9922) -- (3.0610,4.9556) -- (3.1381,4.9062) -- (3.2092,4.8496) -- (3.2822,4.7795) -- (3.3525,4.6992) -- (3.4199,4.6093) -- (3.4804,4.5158) -- (3.5409,4.4082) -- (3.5945,4.2990) -- (3.6926,4.0530) -- (3.7749,3.7749) -- (3.5347,3.8266) -- (3.2762,3.8652) -- (3.0104,3.8884) -- (2.7500,3.8958) -- (2.4806,3.8879) -- (2.2238,3.8652) -- (1.9653,3.8266) -- (1.7251,3.7749) -- cycle;
\fill[even odd rule,blue!11] (1.7251,3.7749) -- (1.6734,3.5347) -- (1.6348,3.2762) -- (1.6116,3.0104) -- (1.6042,2.7500) -- (1.6121,2.4806) -- (1.6348,2.2238) -- (1.6734,1.9653) -- (1.7251,1.7251) -- (1.4470,1.8074) -- (1.2010,1.9055) -- (1.0918,1.9591) -- (0.9842,2.0196) -- (0.8855,2.0838) -- (0.8008,2.1475) -- (0.7247,2.2138) -- (0.6540,2.2867) -- (0.5969,2.3576) -- (0.5468,2.4346) -- (0.5078,2.5132) -- (0.4799,2.5930) -- (0.4634,2.6736) -- (0.4583,2.7500) -- (0.4641,2.8309) -- (0.4799,2.9070) -- (0.5059,2.9824) -- (0.5444,3.0610) -- (0.5938,3.1381) -- (0.6540,3.2133) -- (0.7247,3.2862) -- (0.8008,3.3525) -- (0.8907,3.4199) -- (0.9842,3.4804) -- (1.0856,3.5377) -- (1.2010,3.5945) -- (1.4470,3.6926) -- (1.7251,3.7749) -- cycle;
\fill[even odd rule,blue!9] (3.7749,1.7251) -- (3.6926,1.4470) -- (3.5945,1.2010) -- (3.5409,1.0918) -- (3.4804,0.9842) -- (3.4199,0.8907) -- (3.3525,0.8008) -- (3.2822,0.7205) -- (3.2133,0.6540) -- (3.1381,0.5938) -- (3.0654,0.5468) -- (2.9868,0.5078) -- (2.9070,0.4799) -- (2.8309,0.4641) -- (2.7500,0.4583) -- (2.6691,0.4641) -- (2.5930,0.4799) -- (2.5132,0.5078) -- (2.4390,0.5444) -- (2.3661,0.5908) -- (2.2908,0.6504) -- (2.2178,0.7205) -- (2.1475,0.8008) -- (2.0801,0.8907) -- (2.0196,0.9842) -- (1.9591,1.0918) -- (1.9055,1.2010) -- (1.8074,1.4470) -- (1.7251,1.7251) -- (1.9653,1.6734) -- (2.2238,1.6348) -- (2.4896,1.6116) -- (2.7500,1.6042) -- (3.0194,1.6121) -- (3.2762,1.6348) -- (3.5347,1.6734) -- (3.7749,1.7251) -- cycle;
\draw[black,line width=.65pt] (3.8958,2.7500) -- (3.8879,3.0194) -- (3.8631,3.2937) -- (3.8234,3.5516) -- (3.7689,3.7984) -- (3.7340,3.9243) -- (3.6977,4.0381) -- (3.6581,4.1475) -- (3.6154,4.2520) -- (3.5697,4.3513) -- (3.5178,4.4510) -- (3.4664,4.5385) -- (3.4125,4.6197) -- (3.3564,4.6945) -- (3.2981,4.7625) -- (3.2338,4.8274) -- (3.1718,4.8807) -- (3.1084,4.9267) -- (3.0437,4.9651) -- (2.9780,4.9959) -- (2.9070,5.0201) -- (2.8354,5.0353) -- (2.7680,5.0414) -- (2.7005,5.0395) -- (2.6332,5.0297) -- (2.5663,5.0120) -- (2.4957,4.9845) -- (2.4303,4.9507) -- (2.3661,4.9092) -- (2.2867,4.8460) -- (2.2059,4.7668) -- (2.1285,4.6752) -- (2.0548,4.5717) -- (1.9855,4.4570) -- (1.9240,4.3383) -- (1.8643,4.2038) -- (1.8125,4.0677) -- (1.7637,3.9166) -- (1.7230,3.7663) -- (1.6863,3.6020) -- (1.6575,3.4410) -- (1.6338,3.2675) -- (1.6176,3.0996) -- (1.6074,2.9208) -- (1.6042,2.7500) -- (1.6077,2.5702) -- (1.6176,2.4004) -- (1.6348,2.2238) -- (1.6575,2.0590) -- (1.6880,1.8897) -- (1.7230,1.7337) -- (1.7660,1.5757) -- (1.8125,1.4323) -- (1.8671,1.2892) -- (1.9240,1.1617) -- (1.9889,1.0370) -- (2.0548,0.9283) -- (2.1285,0.8248) -- (2.2059,0.7332) -- (2.2826,0.6577) -- (2.3661,0.5908) -- (2.4303,0.5493) -- (2.4957,0.5155) -- (2.5619,0.4894) -- (2.6287,0.4712) -- (2.6960,0.4609) -- (2.7635,0.4585) -- (2.8309,0.4641) -- (2.9025,0.4787) -- (2.9691,0.5006) -- (3.0350,0.5303) -- (3.1041,0.5705) -- (3.1676,0.6160) -- (3.2297,0.6688) -- (3.2901,0.7289) -- (3.3525,0.8008) -- (3.4089,0.8751) -- (3.4629,0.9559) -- (3.5145,1.0430) -- (3.5634,1.1359) -- (3.6125,1.2413) -- (3.6554,1.3454) -- (3.6952,1.4545) -- (3.7317,1.5680) -- (3.7668,1.6936) -- (3.8219,1.9400) -- (3.8631,2.2063) -- (3.8879,2.4806) -- (3.8958,2.7500);
\draw[black,line width=.65pt] (5.0417,2.7500) -- (5.0391,2.8040) -- (5.0306,2.8623) -- (5.0162,2.9204) -- (4.9976,2.9735) -- (4.9741,3.0262) -- (4.9430,3.0826) -- (4.9062,3.1381) -- (4.8672,3.1885) -- (4.7711,3.2901) -- (4.6555,3.3866) -- (4.5215,3.4769) -- (4.3705,3.5602) -- (4.2038,3.6357) -- (4.0232,3.7027) -- (3.8303,3.7605) -- (3.6270,3.8086) -- (3.4152,3.8465) -- (3.1971,3.8738) -- (2.9746,3.8903) -- (2.7500,3.8958) -- (2.5254,3.8903) -- (2.3029,3.8738) -- (2.0848,3.8465) -- (1.8730,3.8086) -- (1.6697,3.7605) -- (1.4768,3.7027) -- (1.2962,3.6357) -- (1.1295,3.5602) -- (0.9785,3.4769) -- (0.8445,3.3866) -- (0.7289,3.2901) -- (0.6328,3.1885) -- (0.5908,3.1339) -- (0.5570,3.0826) -- (0.5281,3.0306) -- (0.5024,2.9735) -- (0.4838,2.9204) -- (0.4694,2.8623) -- (0.4609,2.8040) -- (0.4583,2.7500) -- (0.4609,2.6960) -- (0.4694,2.6377) -- (0.4825,2.5841) -- (0.5024,2.5265) -- (0.5281,2.4694) -- (0.5570,2.4174) -- (0.5938,2.3619) -- (0.6328,2.3115) -- (0.7289,2.2099) -- (0.8445,2.1134) -- (0.9785,2.0231) -- (1.1295,1.9398) -- (1.2962,1.8643) -- (1.4768,1.7973) -- (1.6697,1.7395) -- (1.8730,1.6914) -- (2.0848,1.6535) -- (2.3029,1.6262) -- (2.5254,1.6097) -- (2.7500,1.6042) -- (2.9746,1.6097) -- (3.1971,1.6262) -- (3.4152,1.6535) -- (3.6270,1.6914) -- (3.8303,1.7395) -- (4.0232,1.7973) -- (4.2038,1.8643) -- (4.3705,1.9398) -- (4.5215,2.0231) -- (4.6555,2.1134) -- (4.7711,2.2099) -- (4.8672,2.3115) -- (4.9092,2.3661) -- (4.9430,2.4174) -- (4.9719,2.4694) -- (4.9976,2.5265) -- (5.0175,2.5841) -- (5.0306,2.6377) -- (5.0391,2.6960) -- (5.0417,2.7500);
\draw[gray!50,thin] (0,0) rectangle (5.5,5.5);
\node[font=\small,fill=white,fill opacity=.85,text opacity=1,inner sep=1pt] at (2.7500,2.7500) {$S_{1}$};
\node[font=\small,fill=white,fill opacity=.85,text opacity=1,inner sep=1pt] at (2.7500,4.4688) {$S_{2}$};
\node[font=\small,fill=white,fill opacity=.85,text opacity=1,inner sep=1pt] at (2.7500,1.0312) {$S_{3}$};
\node[font=\small,fill=white,fill opacity=.85,text opacity=1,inner sep=1pt] at (1.0312,2.7500) {$S_{4}$};
\node[font=\small,fill=white,fill opacity=.85,text opacity=1,inner sep=1pt] at (4.4688,2.7500) {$S_{5}$};
\node[font=\small,fill=white,fill opacity=.85,text opacity=1,inner sep=1pt] at (4.9042,4.9042) {$S_{6}$};
\end{tikzpicture}
\par\smallskip (3)\end{minipage}
\caption{The two-conic arrangements in \cref{ex:even-conics}.
The line at infinity is not a component. Repeated labels in (2)
denote parts of the same projective region.}
\label{fig:even-conics}
\end{figure}
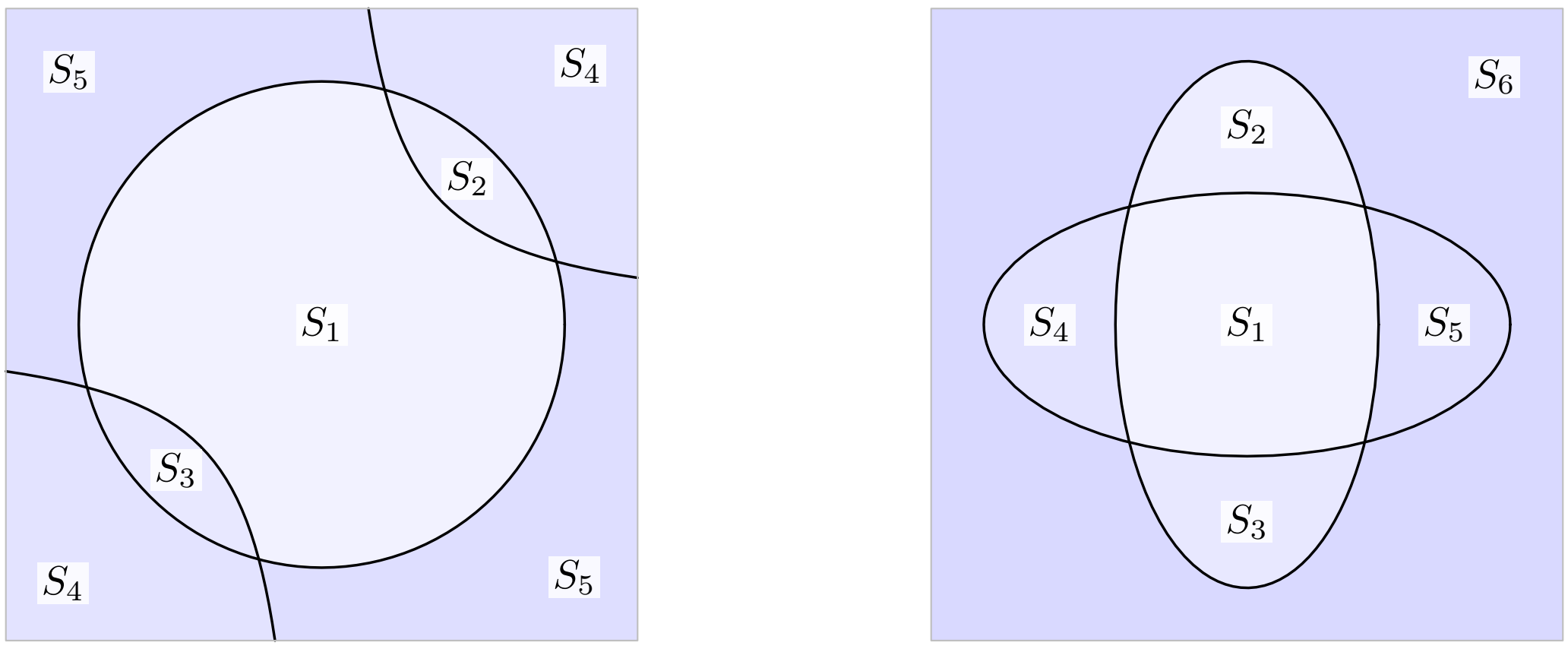

\begin{example}\label{ex:even-conics}
\begin{enumerate}
\item A smooth real conic with a real point has $D=0$ and $\epsilon=0$.
Its two regions are an open disk and an open M\"obius band. The
complexified complement has integral homology
\[
 H_0(\UC;\Z)\simeq\Z,\qquad
 H_1(\UC;\Z)\simeq\Z/2\Z,\qquad H_2(\UC;\Z)=0,
\]
and hence total rational Betti number one.
\item Consider the projective closures of
\[
 xy=1,\qquad x^2+y^2=4.
\]
Their four intersections are real and transverse, and there are no
intersections at infinity. Thus the arrangement is totally real and
$D=4$. An arc of the projective hyperbola joining its two affine
branches across one point at infinity, together with an arc of the
circle, gives a loop in $\CR$ meeting the line at infinity once.
Hence $\epsilon=1$ and $\#\ch(\A)=5$. In \cref{fig:even-conics}(2),
each of $S_4$ and $S_5$ appears in two pieces that join across infinity.
\item Consider the projective closures of
\[
 4x^2+y^2=1,\qquad x^2+4y^2=1.
\]
Their four transverse intersections are
$(\pm1/\sqrt5,\pm1/\sqrt5)$, so $D=4$. The real union lies in the
affine plane; hence $\epsilon=0$ and $\#\ch(\A)=6$.
In \cref{fig:even-conics}(3), $S_6$ is an open M\"obius band and the
other five regions are open disks.
\end{enumerate}
The line at infinity is not a component in these examples. In (2)
and (3), the complexified complements have integral homology
\[
 H_0(\UC;\Z)\simeq\Z,\qquad
 H_1(\UC;\Z)\simeq\Z\oplus\Z/2\Z,\qquad
 H_2(\UC;\Z)\simeq\Z^3,
\]
so their rational Betti numbers are $(1,1,3)$.
These two complexified complements are also homeomorphic. Indeed, ordered
pairs of smooth conics meeting transversely form a connected
Zariski-open subset of $\PP^5_{\C}\times\PP^5_{\C}$. Thom's first
isotopy lemma, applied to the corresponding family of stratified
projective planes, shows that the pairs $(\CP^2,\CC)$ are locally
topologically trivial over this parameter space
\cite[Part I, Section 1.5]{goresky-macpherson}.
\end{example}

\begin{corollary}\label{cor:euler-regions}
For every totally real arrangement of rational curves,
\[
 \sum_{S\in\ch(\A)}\chi(S)=\sum_{k=0}^2 b_k(\UC;\Q)=D+1.
\]
\end{corollary}
\begin{proof}
By \cref{lem:graph-complement}, the left-hand side counts the disk
regions. The assertion follows from \cref{thm:regions}.
\end{proof}

\subsection{Affine curves and bounded regions}
\label{subsec:affine}
Assume $C_0=L_\infty$ is a real line, and let
\[
 C^{\mathrm{aff}}_{\C}=\bigcup_{i=1}^n(C_i(\C)\setminus L_\infty),
 \qquad C^{\mathrm{aff}}_{\R}=C^{\mathrm{aff}}_{\C}\cap\R^2.
\]
The projective arrangement, including $L_\infty$, is assumed totally
real. Define
\begin{equation}\label{eq:affine-data}
 t_i=\#\nu_i^{-1}(C_i\cap L_\infty),\quad
 k=\sum_{i=1}^n t_i,\quad
 D_{\mathrm{aff}}=\sum_{p\notin L_\infty}(r_p-1),\quad
 e=b_0(C^{\mathrm{aff}}_{\R}).
\end{equation}
Every point counted by $t_i$ is real. Thus the real affine normalization
is a disjoint union of $k$ copies of $\R$. A region is bounded if its
projective closure avoids $L_\infty$.

\begin{theorem}\label{thm:bounded}
With the preceding assumptions and $n\geq1$,
\begin{align}
 \#\ch(\A)&=1+k+D_{\mathrm{aff}},\label{eq:aff-total}\\
 \#\bch(\A)&=D_{\mathrm{aff}}-k+e,\label{eq:aff-bound}\\
 \chi(\UC)&=1-2n+k+D_{\mathrm{aff}}.\label{eq:aff-chi}
\end{align}
Consequently
\begin{equation}\label{eq:bounded-correction}
 \#\bch(\A)-\chi(\UC)=e-1-2(k-n).
\end{equation}
In particular, if $t_i=1$ for all $i$ and $C^{\mathrm{aff}}_{\R}$ is
connected, then $\#\bch(\A)=\chi(\UC)$.
\end{theorem}
\begin{proof}
At a point of $L_\infty$, the total branch number is one plus the number
of branches of the affine components. Hence
$D=D_{\mathrm{aff}}+k$. Since $r=n+1$, \cref{thm:regions,prop:betti}
give \eqref{eq:aff-total} and \eqref{eq:aff-chi}.

For the bounded-region formula, first consider the real affine
normalization
\[
 \widetilde C^{\mathrm{aff}}_{\R}
 =\coprod_{i=1}^n
 \left(\widetilde C_i(\R)\setminus
       \nu_i^{-1}(C_i\cap L_\infty)\right)
 =\coprod_{j=1}^k A_j,\qquad A_j\simeq\R.
\]
The circle $\widetilde C_i(\R)$ contributes $t_i$ of these components.
The restricted normalization map
$\nu:\widetilde C^{\mathrm{aff}}_{\R}\to C^{\mathrm{aff}}_{\R}$
is proper. On each $A_j$, the squared distance
$\|\nu(t)\|^2$ tends to infinity at both ends and is monotone sufficiently
far out on either end, since it is semialgebraic.
Choose a closed disk $B_R$ of sufficiently large radius so that it
contains all finite singular points and all bounded complementary
regions, and so that
\[
 I_j=(\nu|_{A_j})^{-1}(B_R)
\]
is a compact interval for every $j$. Enlarge $R$ if necessary to make
the boundary transverse to the curve.

The finite graph $G=C^{\mathrm{aff}}_{\R}\cap B_R$ is obtained from
these $k$ intervals by the following identifications. For each finite
multibranch point $p$, identify the $r_p$ points of $\nu^{-1}(p)$ to
one point. All these preimages lie in the intervals, by the choice of
$R$. This reduces the Euler characteristic by $r_p-1$ and gives
\[
 \chi(G)=k-D_{\mathrm{aff}}.
\]
Outside $B_R$ there are only disjoint ends, each attached to one
boundary point of $G$. Retracting these ends shows that $G$ has the
same homotopy type as $C^{\mathrm{aff}}_{\R}$; in particular,
$b_0(G)=e$. Thus $b_1(G)=D_{\mathrm{aff}}-k+e$.
For a compact graph in the plane, the number of bounded complementary
components is $b_1(G)$, by Alexander duality. By the choice of $B_R$, these are also the bounded components of
$\R^2\setminus C^{\mathrm{aff}}_{\R}$. This proves
\eqref{eq:aff-bound}.
\end{proof}

\begin{remark}\label{rem:essential}
When $t_i=1$ for all $i$, the exact correction is $e-1$. For line
arrangements in the plane with at least two lines, connectedness of the
union is equivalent to the existence of two nonparallel lines, hence to
essentiality. A single line is a harmless exception: both sides of the
bounded formula vanish. For $n$ parallel lines, $e=n$, $\chi(\UC)=1-n$,
and there are no bounded regions.
\end{remark}

\section{Regions and Borel--Moore homology}
\label{sec:BM}
An orientation of a real region $S$ determines a locally finite
fundamental cycle. Since the inclusion $S\hookrightarrow\UC$ is
proper, this cycle defines a class $[S]\in H_2^{\BM}(\UC;\Z)$.
We call the map sending an oriented region to this class the
\emph{region cycle map}.

Let $\ch_{\mathrm{or}}(\A)$ be the set of orientable regions. For
$K=\Z$ or $\Q$, write
\[
 A_{\mathrm{or}}(\A;K)=\bigoplus_{S\in\ch_{\mathrm{or}}(\A)}K e_S.
\]
For each orientable region, choose an orientation. 
By sending $e_S$ to
its fundamental class, one can identify this group with
$H_2(\RP^2,\CR;K)=H_2^{\BM}(\UR;K)$. A nonorientable region has no
fundamental class with these constant coefficients and contributes zero
to this relative group. We obtain
\begin{equation}\label{eq:region-map}
 \gamma_K:A_{\mathrm{or}}(\A;K)\longrightarrow
 H_2(\CP^2,\CC;K)=H_2^{\BM}(\UC;K).
\end{equation}
For a curve arrangement, we define its chamber algebra by
\[
 \VG(\A)_K=\operatorname{Map}(\ch(\A),K)
\]
with pointwise multiplication. If all regions are orientable, its
additive group is identified with $A_{\mathrm{or}}(\A;K)$ by
$\one_S\mapsto e_S$. Reversing the chosen orientation of $S$ replaces
$\gamma_K(e_S)$ by its negative. The map on relative homology is
canonical; this description in terms of chamber functions depends on
the chosen orientations.

\subsection{Rational coefficients}
\begin{theorem}\label{thm:rational-map}
For every totally real arrangement of rational curves, $\gamma_{\Q}$
is surjective. If $\partial_{\R}$ denotes the boundary map for
$(\RP^2,\CR)$, then
\begin{equation}\label{eq:kernelQ}
 \partial_{\R}:\ker\gamma_{\Q}\xrightarrow{\ \simeq\ }
 \bigoplus_{i=0}^n\Q\alpha_i.
\end{equation}
In particular the kernel has dimension $r$.
\end{theorem}
\begin{proof}
The boundary maps fit into the commutative square
\[
\begin{tikzpicture}[baseline=(current bounding box.center),>=Latex]
\node (a) at (0,1.4) {$H_2(\RP^2,\CR;\Q)$};
\node (b) at (6,1.4) {$H_1(\CR;\Q)$};
\node (c) at (0,0) {$H_2(\CP^2,\CC;\Q)$};
\node (d) at (6,0) {$H_1(\CC;\Q)$};
\draw[->] (a)--node[above]{$\partial_{\R}$}node[below]{$\simeq$}(b);
\draw[->] (c)--node[above]{$\partial_{\C}$}node[below]{$\simeq$}(d);
\draw[->] (a)--node[left]{$\gamma_{\Q}$}(c);
\draw[->] (b)--node[right]{$j_*$}(d);
\end{tikzpicture}
\]
The top arrow is an isomorphism since $H_2(\RP^2;\Q)=H_1(\RP^2;\Q)=0$.
The bottom arrow is an isomorphism by \eqref{eq:BM-extension}.
The right arrow is surjective with kernel $L\otimes\Q$, by
\eqref{eq:normalization-sequence}.
\end{proof}

\subsection{Integral coefficients and the degree obstruction}
We refine \eqref{eq:kernelQ} over $\Z$ and determine the obstruction
to integral generation when $\bar d$ is odd.

\begin{theorem}\label{thm:integral}
Suppose $\bar d$ is odd. After choosing orientations of all regions,
let $A=\bigoplus_{S\in\ch(\A)}\Z e_S$. Then
\begin{enumerate}
\item the boundary identifies the kernel of $\gamma_{\Z}$ with
\begin{equation}\label{eq:kernel-lattice}
 L\cap\ker\eta
 =\left\{\sum_{i=0}^n a_i\alpha_i:
                   \sum_{i=0}^n a_i d_i\equiv0\pmod2\right\};
\end{equation}
\item $\im\gamma_{\Z}$ maps isomorphically to $H_1(\CC;\Z)$ under
$\partial_{\C}$;
\item there is an exact sequence
\begin{equation}\label{eq:integer-exact}
 0\longrightarrow L\cap\ker\eta\longrightarrow A
 \xrightarrow{\gamma_{\Z}}H_2^{\BM}(\UC;\Z)
 \longrightarrow\Z/\bar d\Z\longrightarrow0.
\end{equation}
\end{enumerate}
In particular, $\gamma_{\Z}$ is surjective when $\bar d=1$, and its
kernel is then a primitive sublattice of $A$ of rank $r$.
\end{theorem}
\begin{proof}
The real relative sequence gives an isomorphism 
\begin{equation}\label{eq:real-relative-Z}
 A=H_2(\RP^2,\CR;\Z)
 \xrightarrow[\partial_{\R}]{\ \simeq\ }\ker\eta
 \subset H_1(\CR;\Z),
\end{equation}
because $H_2(\RP^2;\Z)=0$. Choose an index $i$ with $d_i$ odd.
By \eqref{eq:normalization-sequence}, every class in $H_1(\CC;\Z)$
has a lift $u\in H_1(\CR;\Z)$. Since $j_*\alpha_i=0$,
adding $\alpha_i$ if necessary makes $\eta(u)=0$ without changing
$j_*u$. 
Thus $\partial_{\C}\gamma_{\Z}$ is surjective.

It remains to control the torsion kernel of $\partial_{\C}$. Complex
conjugation acts as $-1$ on $H_2(\CP^2;\Z)$, hence as $-1$ on its
image $T=\Z/\bar d\Z$ in \eqref{eq:BM-extension}. It fixes the image
of $H_2(\RP^2,\CR;\Z)$ pointwise. Consequently an element of
$\im\gamma_{\Z}\cap T$ is killed by $2$; since $\bar d$ is odd, it
is zero. This proves the second assertion. In the commutative square
of boundary maps, a class maps to zero precisely when its real boundary
belongs to $\ker j_*=L$. Equations \eqref{eq:parity} and
\eqref{eq:real-relative-Z} give \eqref{eq:kernel-lattice}.

Finally, $\im\gamma_{\Z}$ is a section of
\eqref{eq:BM-extension}, so its cokernel is $T$. The parity lattice in
\eqref{eq:kernel-lattice} has rank $r$. When $\bar d=1$, the target is
free, which also proves primitivity of the kernel.
\end{proof}

\begin{example}
For a single real nodal cubic with two real branches at its node,
$r=1$, $D=1$, and $\bar d=3$. Thus
\[
 H_2^{\BM}(\UC;\Z)\simeq\Z\oplus\Z/3\Z,
 \qquad \operatorname{coker}\gamma_{\Z}\simeq\Z/3\Z.
\]
The region classes generate over $\Q$ but not over $\Z$.
\end{example}

\subsection{Winding functions in an affine chart}
\label{subsec:winding}
Suppose $C_0=L_\infty$ and orient $\R^2$ by $dx\wedge dy$. 
Orient each
circle $\widetilde C_i(\R)$ for $1\leq i\leq n$. Away from singular
points this orients the real affine branches of $C_i$.

There exists an integer-valued locally constant function $w_i$ on
$\R^2\setminus C_i(\R)$ whose jump from the right side to the left
side of an oriented branch is $+1$. It is unique up to an additive
constant. Indeed, the oriented real curve is a locally finite
one-cycle in $\R^2$, and $H_1^{\BM}(\R^2;\Z)=0$; a locally finite
two-chain bounding it is constant on each complementary region.
This construction includes self-nodes and cusps by following the
normalization. Restrict $w_i$ to the chambers of the full arrangement.
For an oriented affine line, choose its left side to be positive and
normalize $w_i$ to be zero on the right side. Then $w_i$ is exactly
the Heaviside function of that line.

\begin{corollary}\label{cor:winding}
If the line at infinity is a component, then
\begin{equation}\label{eq:winding-kernel}
 \ker\gamma_{\Z}=\Z\one\oplus\Z w_1\oplus\cdots\oplus\Z w_n.
\end{equation}
Equivalently, all integral relations among region classes are generated
by
\[
 \sum_S[S]=0,\qquad \sum_S w_i(S)[S]=0\quad(1\leq i\leq n).
\]
\end{corollary}
\begin{proof}
Orient the normalization of $C_0$ so that
$\partial_{\R}(\one)=2\alpha_0$. Compactifying the chain for $w_i$
shows that
\[
 \partial_{\R}(w_i)=\alpha_i+m_i\alpha_0
 \quad\text{for some }m_i\in\Z.
\]
The coefficient of $\alpha_i$ is one by construction. The parity
condition forces $m_i\equiv d_i\pmod2$. Thus these $n+1$ boundaries
form a basis of the lattice in \eqref{eq:kernel-lattice}: subtracting
the $w_i$ boundaries removes all $\alpha_i$ terms for $i>0$, and the
remaining coefficient of $\alpha_0$ is even. Apply
\cref{thm:integral}.
\end{proof}
A winding function need not be a $\{0,1\}$-valued polynomial sign
function. Consider the hyperbola $xy=1$ in panel (E1) of
\cref{fig:first}, oriented by increasing $x$ on both affine branches.
Its winding function can be chosen to take the values $-1,0,1$ on
$S_1$ (southwest), $S_2$ (middle), and $S_3$ (northeast), respectively. In particular,
if their classes are $a,b,c$, respectively, then
\[
 a+b+c=0,\qquad -a+c=0.
\]
Hence $c=a$ and $b=-2a$ in $H_2^{\BM}\simeq\Z$. Merely recording the
sign of $xy-1$ would not recover this relation.

\section{Combinatorial and real combinatorial ranks}
\label{sec:positive}
\subsection{Genus zero and canonical forms}
For a complex projective variety $X$ of dimension $m$ and a closed
subvariety $Y$, Brown--Dupont \cite{bd-pos} define
\[
 g(X,Y)=\sum_{p>0}h^{-p,0}(H_m(X,Y;\Q)),\qquad
 \crk(X,Y)=\dim\gr^W_0H_m(X,Y;\Q).
\]
If $X\setminus Y$ is smooth and the genus is zero, their construction
gives a canonical-form map
\[
 \can:H_m(X,Y;\Q)\longrightarrow\Omega^m_{\log}(X\setminus Y),
\]
whose $\C$-linear extension is surjective
\cite[Remark 2.7]{bd-pos}:
\[
 \can_{\C}:H_m(X,Y;\Q)\otimes\C
       \twoheadrightarrow\Omega^m_{\log}(X\setminus Y).
\]
More precisely, $\can$ factors through the weight-zero quotient, and
the induced map
\[
 \gr^W_0H_m(X,Y;\Q)\otimes\C
       \xrightarrow{\ \simeq\ }\Omega^m_{\log}(X\setminus Y)
\]
is an isomorphism. Thus the canonical forms of rational relative
cycles span the space of logarithmic forms over $\C$.
Here logarithmic means logarithmic on a smooth compactification with
normal crossing boundary.

\begin{proposition}\label{prop:pure}
If $\CC\subset\CP^2$ is a nonempty union of geometrically rational
curves, then $H_2(\CP^2,\CC;\Q)$ is pure of type $(0,0)$. In particular,
\begin{equation}\label{eq:genus-rank}
 g(\CP^2,\CC)=0,\qquad \crk(\CP^2,\CC)=D-r+1=h,
\end{equation}
and
\[
 \can_{\C}:H_2(\CP^2,\CC;\Q)\otimes\C
       \xrightarrow{\ \simeq\ }\Omega^2_{\log}(\UC)
\]
is an isomorphism.
\end{proposition}
\begin{proof}
Brown--Dupont \cite[Propositions 3.21--3.22]{bd-pos} show that a union
of rational plane curves has genus zero. Their
Propositions 4.14--4.15 give
$\crk(\CP^2,\CC)=\sum_p(r_p-1)-r+1=h$.
By \cref{prop:betti}, the whole space $H_2(\CP^2,\CC;\Q)$ also has
dimension $h$. Its weight-zero graded piece therefore has the full
dimension, so the relative homology is pure of type $(0,0)$.
The canonical-form isomorphism now follows from their genus-zero
construction.
\end{proof}

In the terminology of Sturmfels--Telen \cite[Definitions 6.2 and
6.10]{st-cub}, a positive arrangement is a genus-zero pair defined over
$\R$, with smooth complement and with smooth real points on the ambient
variety and on every boundary component. Its real combinatorial rank is
the dimension over $\Q$ of the span of the canonical forms of oriented
regions. Our pairs satisfy these conditions. The smooth real points on
the components exist by \cref{def:total}.

\begin{theorem}\label{thm:real-rank}
For a totally real arrangement of rational plane curves,
\begin{equation}\label{eq:rank-equality}
 \crk_{\R}(\CP^2,\CC)=\crk(\CP^2,\CC)=b_2(\UC)=h.
\end{equation}
All rational relations among the canonical forms of orientable regions
are exactly the relations in \eqref{eq:kernelQ}. If the arrangement
contains a line at infinity, they are generated by the constant function
and the winding functions in \cref{cor:winding}.
\end{theorem}
\begin{proof}
Compose the surjection of \cref{thm:rational-map} with the injective
canonical-form map in \cref{prop:pure}. Its image has dimension $h$
over $\Q$ and its complex span is all of $\Omega^2_{\log}(\UC)$.
Injectivity of the latter map also proves the assertion about relations.
\end{proof}
Thus the degree condition needed to count \emph{all} real chambers is
not needed for equality of the two combinatorial ranks. When an open
M\"obius band occurs, it is excluded from the oriented-region cycle
space, exactly as in the definition of regions for positive arrangements.

For hyperplane arrangements, the chamber classes and their differential
form representatives were studied in \cite{yos-ch}. Eur--Lam
\cite{eur-lam} develop canonical forms in the more general setting of
oriented matroids. See also \cite{lam-mat,lam-icm} for their relation
to scattering amplitudes. For curves, residues on the normalizations
and the adjoint-curve methods of \cite{tel-poly} can be used to compute
individual canonical forms. The results above determine the dimension
of their span and all linear relations among them.

\section{Examples and explicit relations}
\label{sec:examples}
The arrangements (E1)--(E8) are totally real and include the line
at infinity $C_0=\{z=0\}$. We give their equations in the affine
chart $z=1$.

\subsection{Conics, a nodal cubic, and a cusp}
Consider the following arrangements of affine curves:
\begin{align*}
 \mathrm{(E1)}&\quad xy=1;\\
 \mathrm{(E2)}&\quad y=0,\quad y=x^2;\\
 \mathrm{(E3)}&\quad y=x^2,\quad y=(x-1)^2,\quad y=(x+1)^2;\\
 \mathrm{(E4)}&\quad xy=1,\quad (x+1)(y-1)=1;\\
 \mathrm{(E5)}&\quad y^2=x^2(x+1);\\
 \mathrm{(E6)}&\quad y=x^2,\quad y^2=x^3.
\end{align*}
Their numerical data are collected in \cref{tab:data}.
The last two rows refer to the Chebyshev curves discussed below.

\begin{table}[htbp]
\centering\small
\begin{tabular}{lrrrrrrrr}
\toprule
 & $n$ & $k$ & $D_{\rm aff}$ & $e$ & $b_2$ & $\#\ch$ & $\#\bch$ & $\chi(\UC)$\\
\midrule
(E1) hyperbola &1&2&0&2&1&3&0&1\\
(E2) line and tangent parabola &2&2&1&1&1&4&0&0\\
(E3) three parabolas &3&3&3&1&3&7&1&1\\
(E4) two hyperbolas &2&4&2&2&4&7&0&3\\
(E5) nodal cubic &1&1&1&1&1&3&1&1\\
(E6) parabola and cusp &2&2&2&1&2&5&1&1\\
(E7) Chebyshev $(3,4)$ &1&1&3&1&3&5&3&3\\
(E8) Chebyshev $(4,5)$ &1&1&6&1&6&8&6&6\\
\bottomrule
\end{tabular}
\caption{Data for the arrangements including the line at infinity.
In every row $b_0=1$, $b_1=n$, and $b_2=k+D_{\rm aff}-n$.}
\label{tab:data}
\end{table}

For (E2), the tangency at $(0,0)$ contributes $1$ to $D_{\rm aff}$,
independently of its intersection multiplicity $2$. For (E3), the three
finite intersections are
$(-1/2,1/4)$, $(1/2,1/4)$, and $(0,1)$; all three parabolas meet the
line at infinity at the same point, producing four local branches in
the projective union.

For (E4), eliminating $y$ gives $x^2+x-1=0$, so both finite
intersections are real and transverse. At each of
$[1:0:0]$ and $[0:1:0]$, the two conics and the line at infinity give
three branches. The two affine connected components each contain one
branch of each hyperbola. Thus $e=2$ and there are no bounded regions.

For (E5), the normalization is
$(x,y)=(t^2-1,t(t^2-1))$. The point $(0,0)$ has the two real preimages
$t=\pm1$, and the point at infinity has one preimage. For (E6), the
cusp has normalization $(x,y)=(t^2,t^3)$. The two finite intersections
with the parabola are $(0,0)$ and $(1,1)$; each contributes one to
$D_{\rm aff}$. The cusp at the origin is unibranch on the cubic, but
there are two branches in the union with the parabola.

\begin{example}[Failure of total reality]\label{ex:nonreal}
The arrangement $z(x^2+y^2-z^2)=0$ is not totally real: its two
intersection points at infinity are $[1:\pm i:0]$. Its real affine
complement has two regions, whereas its rational Betti numbers are
$(1,1,1)$. The two real regions are a disk and the exterior of a
circle in $\R^2$, with Euler characteristics $1$ and $0$. Hence
\cref{cor:euler-regions} also fails without total reality. The cubic $y^2z=x^3-x^2z$ also fails the definition: the
origin in $z=1$ is a solitary real node, with conjugate nonreal tangent
branches. The normalization does not map its real locus onto that
isolated real point. These examples explain why reality of the
singular point alone is not the required condition.
\end{example}

\subsection{Chebyshev curves}
Write $T_m$ for the Chebyshev polynomial, so that
$T_m(\cos\theta)=\cos(m\theta)$. Fix coprime integers $2\leq a<b$.
The Chebyshev curve is
\begin{equation}\label{eq:cheb}
 C_{a,b}:\quad T_b(x)-T_a(y)=0,
 \qquad (x,y)=(T_a(t),T_b(t)).
\end{equation}
These curves are treated in \cite{kp} and in the broader study of
expressive curves \cite{fs}. We include the specific argument needed
for total reality.

\begin{proposition}\label{prop:cheb}
The projective closure of $C_{a,b}$ is an irreducible rational curve of
degree $b$. Its affine singularities are precisely
$(a-1)(b-1)/2$ ordinary nodes, all with two real branches. Its
normalization has one point over infinity. Consequently its projective
closure together with $L_\infty$ is totally real, and, writing
$\delta=(a-1)(b-1)/2$,
\[
 (b_0,b_1,b_2)=(1,1,\delta),\qquad
 \#\ch=\delta+2,\qquad\#\bch=\delta.
\]
\end{proposition}
\begin{proof}
The parametrization extends to a morphism $\phi:\PP^1\to\PP^2$:
writing $T_m^h(u,v)=v^mT_m(u/v)$, it is given by
\[
 [u:v]\longmapsto
 [v^{b-a}T_a^h(u,v):T_b^h(u,v):v^b].
\]
These three coordinates have no common zero. Its image $B$ is a
closed irreducible projective curve. Let $m$ be the number of inverse
images of a general point of $B$. For a general $c\in\C$, the equation
$T_a(t)=c$ has $a$ distinct roots. They occur in groups of $m$, one
for each point of $B$ on the vertical line $x=c$. Thus $m$ divides
$a$. Applying the same argument to $T_b(t)=c$ shows that $m$ divides
$b$. Since $a$ and $b$ are coprime, $m=1$.

Choose a general line $L:px+qy+r=0$ with $q\ne0$, avoiding the
singular points of $B$ and its point at infinity, and meeting $B$
transversely at smooth points. Over the smooth locus the birational
map $\phi$ is one-to-one. Hence the degree-$b$ polynomial
$pT_a(t)+qT_b(t)+r$ has $b$ simple roots giving $b$ distinct points of
$L\cap B$. All intersection multiplicities are one, so $\deg B=b$. The polynomial
$T_b(x)-T_a(y)$ has degree $b$ and vanishes on $B$, so it defines
precisely this irreducible curve. The morphism $\phi$ is therefore the normalization map of $B$.

Every affine singularity satisfies $T_b'(x)=T_a'(y)=0$. The complete
lists of roots of these derivatives are
\[
 x=\cos\frac{j\pi}{b}\quad(1\leq j<b),\qquad
 y=\cos\frac{l\pi}{a}\quad(1\leq l<a).
\]
The point lies on the curve precisely when $j\equiv l\pmod2$.
There are $(a-1)(b-1)/2$ such pairs. At each one, the two diagonal
entries of the Hessian of $T_b(x)-T_a(y)$ are nonzero and have opposite
signs, because
\[
 T_b''\!\left(\cos\frac{j\pi}{b}\right)
   =\frac{(-1)^{j+1}b^2}{\sin^2(j\pi/b)}.
\]
The real quadratic tangent cone therefore splits into two distinct
real lines. All these singularities are ordinary real nodes and the
critical-point calculation excludes any further affine singularities,
including nonreal ones.

The only pole of the parametrization is $t=\infty$, so the
normalization has exactly one point over infinity and the curve is
unibranch there. The union with $L_\infty$ has two real branches at
that point. Its real affine locus is the image of the connected real
line. Now apply \cref{prop:betti,thm:regions,thm:bounded}.
\end{proof}

The cases $(a,b)=(3,4)$ and $(4,5)$ are (E7) and (E8), with $3$ and
$6$ affine nodes, respectively.

\subsection{Relations for the drawn regions}
\label{subsec:example-relations}
The labels in \cref{fig:first,fig:second,fig:third,fig:fourth} are local
to each panel. In \cref{tab:relations}, a row $(a_1,\ldots,a_m)$ means
the relation $\sum_{j=1}^m a_j[S_j]=0$. In addition to the listed rows,
every panel has the relation $(1,\ldots,1)$. These rows and the
constant row form an integral basis of the kernel by
\cref{cor:winding}.

The winding functions follow the component order in (E1)--(E6).
For graphs $y=f(x)$ the orientation is increasing $x$; for the
hyperbolas it is increasing $x$ on both branches; for (E5), the cusp,
and the Chebyshev curves it is increasing normalization parameter $t$.
An additive constant was chosen separately for each row. The relations
can be verified directly by taking the left-minus-right jump across
each oriented branch in the figures.

\begin{table}[htbp]\centering\small
\begin{tabular}{ll}\toprule
Example & Nonconstant winding rows\\\midrule
(E1) & $(0,1,2)$\\
(E2) & $(0,1,1,1)$\\
 & $(0,0,1,0)$\\
(E3) & $(0,1,0,0,1,1,1)$\\
 & $(0,0,0,1,0,1,1)$\\
 & $(0,0,1,0,1,0,1)$\\
(E4) & $(0,0,1,1,2,2,1)$\\
 & $(0,-1,-1,0,0,1,1)$\\
(E5) & $(0,-1,1)$\\
(E6) & $(0,1,0,1,0)$\\
 & $(0,-1,-1,0,0)$\\
(E7) & $(0,0,1,-1,0)$\\
(E8) & $(0,0,-1,1,-1,-2,0,0)$\\
\bottomrule\end{tabular}
\caption{Integral generators for region relations, together with the constant row.}\label{tab:relations}
\end{table}

\subsection{Six points, fifteen lines, and six conics}
\label{subsec:sixpoints}
Take six distinct points $p_1,\ldots,p_6\in\RP^2$. Assume that
\begin{enumerate}
\item[(G1)] no three of the points are collinear;
\item[(G2)] no conic passes through all six points.
\end{enumerate}
Let $L_{ij}$ be the line through $p_i,p_j$, and let $Q_j$ be the
unique conic through the five points other than $p_j$. Condition~(G1)
ensures that each $Q_j$ is smooth, and (G2) ensures that
$p_j\notin Q_j$. We impose two further conditions:
\begin{enumerate}
\item[(G3)] no three distinct lines among the $L_{ij}$ are concurrent
outside $\{p_1,\ldots,p_6\}$;
\item[(G4)] for every $i\ne j$, $L_{ij}$ is not tangent to $Q_j$ at $p_i$.
\end{enumerate}
These conditions define a nonempty Zariski-open subset of the space
of six-point configurations. The arrangement consists of the fifteen
lines $L_{ij}$ and the six conics $Q_i$.

We check the intersections explicitly. Two conics share four of the
base points, so B\'ezout's theorem implies that they meet transversely
at precisely these points. If $k\notin\{i,j\}$, the line $L_{ij}$
meets $Q_k$ transversely at $p_i$ and $p_j$. The line $L_{ij}$ meets
$Q_j$ at $p_i$ and at one further point; (G4) makes the second point
distinct from $p_i$. It is real because the curves are defined over
$\R$ and their other intersection is real. The same applies to $Q_i$.
Thus the arrangement is totally real, and all ten branches at each
base point have distinct tangent directions.

Every nonbase intersection is an ordinary double point. Indeed,
two conics meet only at base points, and (G3) excludes concurrence of
three lines elsewhere. A nonbase point of $Q_k$ can lie on $L_{ij}$
only if $k\in\{i,j\}$; all such lines pass through $p_k\notin Q_k$,
so two of them cannot meet on $Q_k$.

At each $p_i$ there are five lines and five conics, giving a
contribution $9$. Outside the six base points there are
\[
 \binom{15}{2}-6\binom52=45
\]
line--line intersections and $15\cdot2=30$ residual line--conic
intersections, all ordinary double points. Hence
\[
 r=21,\qquad D=6\cdot9+45+30=129.
\]
Our formulas yield
\[
 (b_0,b_1,b_2)=(1,20,109),\quad \#\ch=130,\quad
 \crk_{\R}=\crk=109.
\]
The blow-up at the six points induces an isomorphism from the
complement of the $27$ lines on the associated split smooth cubic
surface to the plane-curve complement. Six of these lines are the
exceptional curves. The classical $130$ regions
are discussed by Klein \cite{klein}; a modern treatment, including
canonical forms, is given in \cite{st-cub}. The above calculation uses
only branch incidence and provides the region-cycle interpretation of
the number $109$.

\section{Further questions}
\label{sec:questions}
In the affine case, the results already give a multiplicative filtration
of the chamber algebra $A=\VG(\A)_{\Q}$ by setting
\[
 \Fil^0A=\Q\one,\qquad \Fil^1A=\ker\gamma_{\Q},\qquad
 \Fil^2A=A.
\]
Its graded dimensions are $(1,n,b_2)$. For line arrangements, the
winding functions are Heaviside functions, so \cref{prop:VG-top,cor:winding}
show that this is the usual VG filtration. For general
curves, multiplicativity follows simply from the fact that the
filtration has length two. In particular, generation in degree one can fail:
in (E7) the sole nonconstant winding function takes the same value on
several distinct chambers, so it and $\one$ cannot generate the full
algebra of chamber functions. It would be interesting to find a geometric description of this
filtration and to determine how it depends on the affine chart.

A second problem is to compute canonical forms for singular boundary
curves directly from normalization and branch incidence. Such forms must satisfy the linear
relations proved here, in addition to the local residue conditions.
Finally, in higher dimensions the finite branch graph is replaced by
more complicated strata and dual complexes. Both the genus-zero
condition and generation by real relative cycles require additional
hypotheses; neither follows formally from the plane-curve argument.

\subsection*{Acknowledgements}
The author thanks Simon Telen for suggesting a comparison between the
combinatorial and real combinatorial ranks, and Sakumi Sugawara for
clarifying the hypotheses needed in the bounded-region formula.

This paper was prepared with substantial assistance from the AI tool
ChatGPT, which provided mathematical suggestions and helped simplify
arguments and clarify the exposition.  The author is
responsible for all assertions in this paper.

This work was partially supported by JSPS KAKENHI Grant Numbers
JP23H00081 and JP25K24692.

\appendix
\section{Affine pictures}
\label{app:pictures}
All arrangements include the line at infinity, which is not drawn.
The frames mark the viewing windows.

\begin{figure}[htbp]\centering
\begin{minipage}{.47\textwidth}\centering
\begin{tikzpicture}[x=1cm,y=1cm,line cap=round,line join=round]
\fill[blue!3] (2.4444,0.0000) -- (0.0000,0.0000) -- (0.0000,2.4444) -- (0.3544,2.3992) -- (0.6647,2.3471) -- (0.9293,2.2885) -- (1.1666,2.2193) -- (1.3673,2.1423) -- (1.5407,2.0552) -- (1.6228,2.0045) -- (1.6958,1.9529) -- (1.7688,1.8936) -- (1.8333,1.8333) -- (1.8966,1.7654) -- (1.9513,1.6979) -- (2.0061,1.6205) -- (2.0608,1.5307) -- (2.1065,1.4443) -- (2.1521,1.3446) -- (2.2251,1.1492) -- (2.2890,0.9273) -- (2.3529,0.6342) -- (2.3985,0.3595) -- (2.4444,0.0000) -- cycle;
\fill[blue!6] (5.5000,3.0556) -- (5.5000,0.0000) -- (2.4444,0.0000) -- (2.3985,0.3595) -- (2.3529,0.6342) -- (2.2890,0.9273) -- (2.2251,1.1492) -- (2.1521,1.3446) -- (2.1065,1.4443) -- (2.0608,1.5307) -- (2.0061,1.6205) -- (1.9513,1.6979) -- (1.8966,1.7654) -- (1.8333,1.8333) -- (1.7688,1.8936) -- (1.6958,1.9529) -- (1.6228,2.0045) -- (1.5407,2.0552) -- (1.3673,2.1423) -- (1.1666,2.2193) -- (0.9293,2.2885) -- (0.6647,2.3471) -- (0.3544,2.3992) -- (0.0000,2.4444) -- (0.0000,5.5000) -- (3.0556,5.5000) -- (3.1015,5.1405) -- (3.1471,4.8658) -- (3.2110,4.5727) -- (3.2749,4.3508) -- (3.3479,4.1554) -- (3.3935,4.0557) -- (3.4392,3.9693) -- (3.4939,3.8795) -- (3.5487,3.8021) -- (3.6034,3.7346) -- (3.6667,3.6667) -- (3.7312,3.6064) -- (3.8042,3.5471) -- (3.8772,3.4955) -- (3.9593,3.4448) -- (4.1327,3.3577) -- (4.3334,3.2807) -- (4.5707,3.2115) -- (4.8353,3.1529) -- (5.1456,3.1008) -- (5.5000,3.0556) -- cycle;
\fill[blue!9] (3.0556,5.5000) -- (5.5000,5.5000) -- (5.5000,3.0556) -- (5.1456,3.1008) -- (4.8353,3.1529) -- (4.5707,3.2115) -- (4.3334,3.2807) -- (4.1327,3.3577) -- (3.9593,3.4448) -- (3.8772,3.4955) -- (3.8042,3.5471) -- (3.7312,3.6064) -- (3.6667,3.6667) -- (3.6034,3.7346) -- (3.5487,3.8021) -- (3.4939,3.8795) -- (3.4392,3.9693) -- (3.3935,4.0557) -- (3.3479,4.1554) -- (3.2749,4.3508) -- (3.2110,4.5727) -- (3.1471,4.8658) -- (3.1015,5.1405) -- (3.0556,5.5000) -- cycle;
\draw[black,line width=.65pt] (0.0000,2.4444) -- (0.3544,2.3992) -- (0.6647,2.3471) -- (0.8015,2.3187) -- (0.9293,2.2885) -- (1.0479,2.2563) -- (1.1666,2.2193) -- (1.2669,2.1834) -- (1.3673,2.1423) -- (1.4586,2.0993) -- (1.5407,2.0552) -- (1.6228,2.0045) -- (1.6958,1.9529) -- (1.7688,1.8936) -- (1.8333,1.8333) -- (1.8966,1.7654) -- (1.9513,1.6979) -- (2.0061,1.6205) -- (2.0608,1.5307) -- (2.1065,1.4443) -- (2.1521,1.3446) -- (2.1886,1.2532) -- (2.2251,1.1492) -- (2.2616,1.0295) -- (2.2890,0.9273) -- (2.3255,0.7706) -- (2.3529,0.6342) -- (2.3985,0.3595) -- (2.4444,0.0000);
\draw[black,line width=.65pt] (3.0556,5.5000) -- (3.1015,5.1405) -- (3.1471,4.8658) -- (3.1745,4.7294) -- (3.2110,4.5727) -- (3.2384,4.4705) -- (3.2749,4.3508) -- (3.3114,4.2468) -- (3.3479,4.1554) -- (3.3935,4.0557) -- (3.4392,3.9693) -- (3.4939,3.8795) -- (3.5487,3.8021) -- (3.6034,3.7346) -- (3.6667,3.6667) -- (3.7312,3.6064) -- (3.8042,3.5471) -- (3.8772,3.4955) -- (3.9593,3.4448) -- (4.0414,3.4007) -- (4.1327,3.3577) -- (4.2331,3.3166) -- (4.3334,3.2807) -- (4.4521,3.2437) -- (4.5707,3.2115) -- (4.6985,3.1813) -- (4.8353,3.1529) -- (5.1456,3.1008) -- (5.5000,3.0556);
\draw[gray!50,thin] (0.0000,0.0000) rectangle (5.5000,5.5000);
\node[font=\small,fill=white,fill opacity=.85,text opacity=1,inner sep=1pt] at (1.0742,1.0742) {$S_{1}$};
\node[font=\small,fill=white,fill opacity=.85,text opacity=1,inner sep=1pt] at (1.6812,3.8188) {$S_{2}$};
\node[font=\small,fill=white,fill opacity=.85,text opacity=1,inner sep=1pt] at (4.4258,4.4258) {$S_{3}$};
\end{tikzpicture}
\par\smallskip (E1)\end{minipage}\hfill
\begin{minipage}{.47\textwidth}\centering
\begin{tikzpicture}[x=1cm,y=1cm,line cap=round,line join=round]
\fill[blue!3] (5.5000,0.8967) -- (5.5000,0.0000) -- (0.0000,0.0000) -- (0.0000,0.8967) -- (5.5000,0.8967) -- cycle;
\fill[blue!6] (0.0000,0.8967) -- (0.0000,4.8424) -- (0.5780,4.8424) -- (0.7222,4.3359) -- (0.8561,3.8967) -- (0.9900,3.4875) -- (1.1239,3.1082) -- (1.2578,2.7590) -- (1.3917,2.4397) -- (1.5257,2.1505) -- (1.6596,1.8912) -- (1.7935,1.6620) -- (1.9274,1.4627) -- (2.0804,1.2717) -- (2.2143,1.1367) -- (2.3483,1.0317) -- (2.4822,0.9567) -- (2.6161,0.9117) -- (2.6926,0.8995) -- (2.7500,0.8967) -- (0.0000,0.8967) -- cycle;
\fill[blue!9] (0.5780,4.8424) -- (4.9220,4.8424) -- (4.7013,4.0813) -- (4.4717,3.3760) -- (4.3570,3.0565) -- (4.2422,2.7590) -- (4.1465,2.5279) -- (4.0317,2.2708) -- (3.9361,2.0733) -- (3.8213,1.8566) -- (3.7065,1.6620) -- (3.5917,1.4893) -- (3.4961,1.3623) -- (3.3813,1.2301) -- (3.2665,1.1199) -- (3.1517,1.0317) -- (3.0752,0.9852) -- (2.9987,0.9485) -- (2.9222,0.9215) -- (2.8457,0.9044) -- (2.7691,0.8970) -- (2.6926,0.8995) -- (2.6161,0.9117) -- (2.5204,0.9408) -- (2.4439,0.9751) -- (2.3674,1.0192) -- (2.2909,1.0730) -- (2.1952,1.1542) -- (2.0422,1.3158) -- (1.9465,1.4367) -- (1.8700,1.5444) -- (1.7170,1.7893) -- (1.5448,2.1116) -- (1.3726,2.4835) -- (1.2196,2.8557) -- (1.0665,3.2671) -- (0.8943,3.7767) -- (0.7413,4.2713) -- (0.5780,4.8424) -- cycle;
\fill[blue!3] (4.9220,4.8424) -- (5.5000,4.8424) -- (5.5000,0.8967) -- (2.7500,0.8967) -- (2.8074,0.8995) -- (2.8839,0.9117) -- (3.0178,0.9567) -- (3.1517,1.0317) -- (3.2857,1.1367) -- (3.4387,1.2934) -- (3.5726,1.4627) -- (3.7065,1.6620) -- (3.8404,1.8912) -- (3.9743,2.1505) -- (4.1083,2.4397) -- (4.2422,2.7590) -- (4.3761,3.1082) -- (4.5100,3.4875) -- (4.6439,3.8967) -- (4.7778,4.3359) -- (4.9220,4.8424) -- cycle;
\draw[black,line width=.65pt] (0.0000,0.8967) -- (5.5000,0.8967);
\draw[blue!65!black,line width=.65pt] (0.5780,4.8424) -- (0.7222,4.3359) -- (0.8561,3.8967) -- (0.9900,3.4875) -- (1.1239,3.1082) -- (1.2578,2.7590) -- (1.3917,2.4397) -- (1.5257,2.1505) -- (1.6596,1.8912) -- (1.7935,1.6620) -- (1.9274,1.4627) -- (2.0039,1.3623) -- (2.0804,1.2717) -- (2.1378,1.2102) -- (2.2143,1.1367) -- (2.2909,1.0730) -- (2.3483,1.0317) -- (2.4248,0.9852) -- (2.4822,0.9567) -- (2.5587,0.9273) -- (2.6161,0.9117) -- (2.6926,0.8995) -- (2.7500,0.8967) -- (2.8074,0.8995) -- (2.8839,0.9117) -- (2.9604,0.9338) -- (3.0178,0.9567) -- (3.0752,0.9852) -- (3.1517,1.0317) -- (3.2091,1.0730) -- (3.2857,1.1367) -- (3.3622,1.2102) -- (3.4387,1.2934) -- (3.5726,1.4627) -- (3.7065,1.6620) -- (3.8404,1.8912) -- (3.9743,2.1505) -- (4.1083,2.4397) -- (4.2422,2.7590) -- (4.3761,3.1082) -- (4.5100,3.4875) -- (4.6439,3.8967) -- (4.7778,4.3359) -- (4.9220,4.8424);
\draw[gray!50,thin] (0.0000,0.0000) rectangle (5.5000,4.8424);
\node[font=\small,fill=white,fill opacity=.85,text opacity=1,inner sep=1pt] at (2.7500,0.4484) {$S_{1}$};
\node[font=\small,fill=white,fill opacity=.85,text opacity=1,inner sep=1pt] at (0.8221,1.7194) {$S_{2}$};
\node[font=\small,fill=white,fill opacity=.85,text opacity=1,inner sep=1pt] at (2.7479,3.2704) {$S_{3}$};
\node[font=\small,fill=white,fill opacity=.85,text opacity=1,inner sep=1pt] at (4.6779,1.7194) {$S_{4}$};
\end{tikzpicture}
\par\smallskip (E2)\end{minipage}
\caption{A hyperbola (E1) and a line tangent to a parabola (E2).}\label{fig:first}
\end{figure}
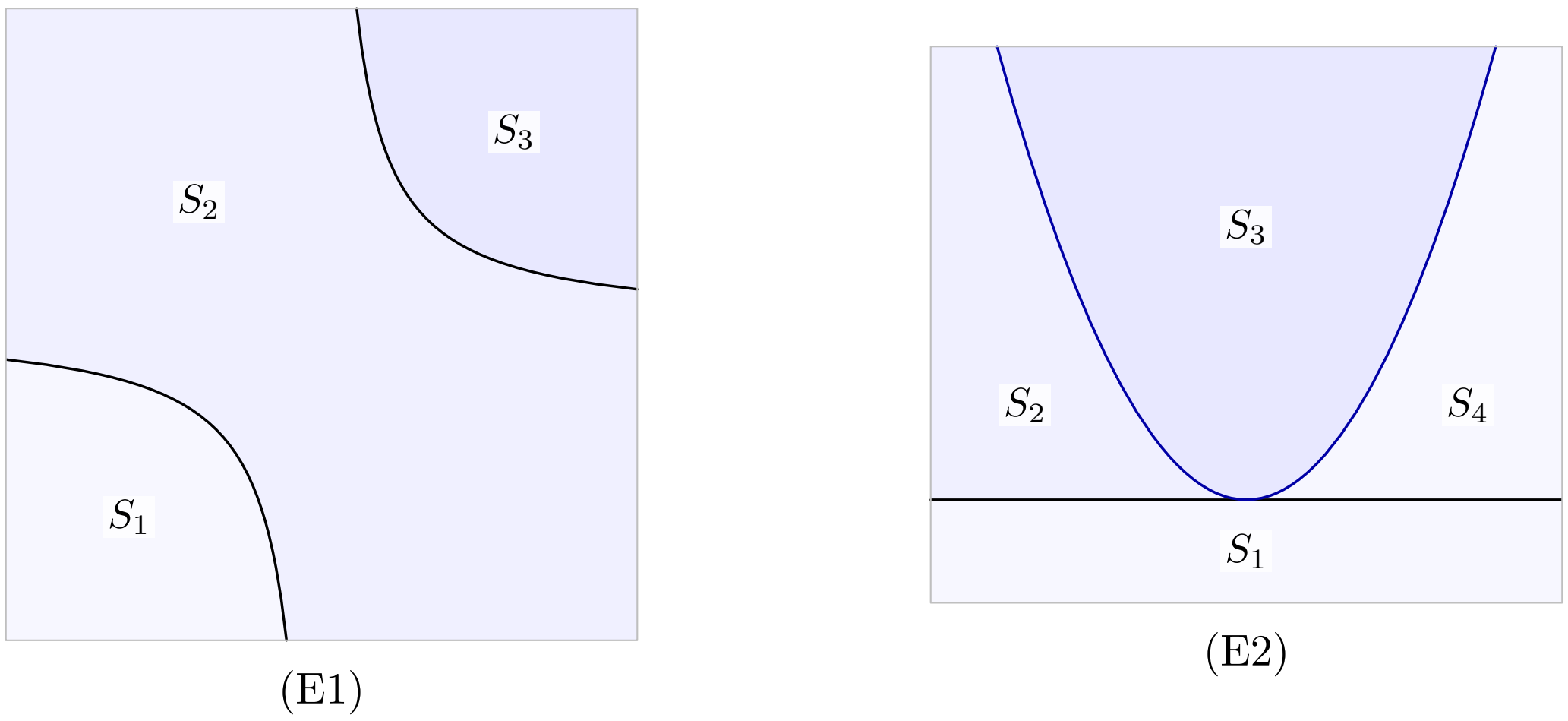
\begin{figure}[htbp]\centering
\begin{minipage}{.47\textwidth}\centering
\begin{tikzpicture}[x=1cm,y=1cm,line cap=round,line join=round]
\fill[blue!3] (5.5000,2.1084) -- (5.5000,0.0000) -- (0.0000,0.0000) -- (0.0000,2.1084) -- (0.1729,1.7509) -- (0.3614,1.4129) -- (0.5500,1.1293) -- (0.7176,0.9228) -- (0.8852,0.7592) -- (0.9690,0.6935) -- (1.0738,0.6265) -- (1.1786,0.5762) -- (1.2624,0.5480) -- (1.3462,0.5306) -- (1.4405,0.5238) -- (1.5138,0.5279) -- (1.5976,0.5427) -- (1.6814,0.5681) -- (1.7652,0.6044) -- (1.8490,0.6513) -- (1.9329,0.7089) -- (2.0952,0.8513) -- (2.2262,0.7333) -- (2.3729,0.6324) -- (2.5195,0.5644) -- (2.6662,0.5292) -- (2.7710,0.5241) -- (2.8548,0.5322) -- (2.9386,0.5510) -- (3.0433,0.5895) -- (3.1271,0.6324) -- (3.2319,0.7012) -- (3.3157,0.7682) -- (3.4048,0.8513) -- (3.5252,0.7418) -- (3.6300,0.6647) -- (3.7348,0.6044) -- (3.8605,0.5541) -- (3.9862,0.5279) -- (4.0910,0.5246) -- (4.1957,0.5380) -- (4.3214,0.5762) -- (4.4681,0.6513) -- (4.6148,0.7592) -- (4.7614,0.9000) -- (4.9081,1.0737) -- (5.0548,1.2802) -- (5.2014,1.5196) -- (5.3481,1.7918) -- (5.5000,2.1084) -- cycle;
\fill[blue!6] (3.4048,0.8513) -- (3.2529,0.7169) -- (3.1690,0.6579) -- (3.0852,0.6096) -- (3.0014,0.5721) -- (2.9176,0.5453) -- (2.8338,0.5292) -- (2.7500,0.5238) -- (2.6662,0.5292) -- (2.5824,0.5453) -- (2.4986,0.5721) -- (2.4148,0.6096) -- (2.3310,0.6579) -- (2.2471,0.7169) -- (2.0952,0.8513) -- (2.2471,1.0207) -- (2.4148,1.2487) -- (2.5824,1.5196) -- (2.7500,1.8333) -- (2.9176,1.5196) -- (3.0852,1.2487) -- (3.2529,1.0207) -- (3.4048,0.8513) -- cycle;
\fill[blue!9] (0.0000,2.1084) -- (0.0000,4.5833) -- (0.4443,4.5833) -- (0.6338,3.9436) -- (0.8433,3.2999) -- (1.0529,2.7233) -- (1.2624,2.2137) -- (1.4719,1.7712) -- (1.5767,1.5751) -- (1.6814,1.3958) -- (1.7862,1.2332) -- (1.8910,1.0873) -- (1.9957,0.9583) -- (2.0952,0.8513) -- (1.9538,0.7250) -- (1.8281,0.6385) -- (1.7024,0.5762) -- (1.5767,0.5380) -- (1.4510,0.5239) -- (1.3252,0.5340) -- (1.1995,0.5681) -- (1.0529,0.6385) -- (0.9271,0.7250) -- (0.8014,0.8357) -- (0.6757,0.9704) -- (0.5290,1.1582) -- (0.4033,1.3452) -- (0.2567,1.5940) -- (0.1310,1.8333) -- (0.0000,2.1084) -- cycle;
\fill[blue!3] (5.0557,4.5833) -- (5.5000,4.5833) -- (5.5000,2.1084) -- (5.3062,1.7106) -- (5.2014,1.5196) -- (5.0967,1.3452) -- (4.9919,1.1877) -- (4.8871,1.0469) -- (4.7824,0.9228) -- (4.6776,0.8156) -- (4.5729,0.7250) -- (4.4681,0.6513) -- (4.3633,0.5943) -- (4.2586,0.5541) -- (4.1538,0.5306) -- (4.0595,0.5238) -- (3.9652,0.5306) -- (3.8605,0.5541) -- (3.7348,0.6044) -- (3.6300,0.6647) -- (3.5252,0.7418) -- (3.4048,0.8513) -- (3.5043,0.9583) -- (3.6090,1.0873) -- (3.7138,1.2332) -- (3.8186,1.3958) -- (3.9233,1.5751) -- (4.0281,1.7712) -- (4.2376,2.2137) -- (4.4471,2.7233) -- (4.6567,3.2999) -- (4.8662,3.9436) -- (5.0557,4.5833) -- cycle;
\fill[blue!6] (0.4443,4.5833) -- (1.7539,4.5833) -- (1.9957,3.7764) -- (2.2471,3.0321) -- (2.3729,2.6962) -- (2.4986,2.3845) -- (2.6243,2.0968) -- (2.7500,1.8333) -- (2.5824,1.5196) -- (2.4148,1.2487) -- (2.2471,1.0207) -- (2.0952,0.8513) -- (1.9957,0.9583) -- (1.8910,1.0873) -- (1.7862,1.2332) -- (1.6814,1.3958) -- (1.5767,1.5751) -- (1.4719,1.7712) -- (1.2624,2.2137) -- (1.0529,2.7233) -- (0.8433,3.2999) -- (0.6338,3.9436) -- (0.4443,4.5833) -- cycle;
\fill[blue!9] (3.7461,4.5833) -- (5.0557,4.5833) -- (4.8662,3.9436) -- (4.6567,3.2999) -- (4.4471,2.7233) -- (4.2376,2.2137) -- (4.0281,1.7712) -- (3.9233,1.5751) -- (3.8186,1.3958) -- (3.7138,1.2332) -- (3.6090,1.0873) -- (3.5043,0.9583) -- (3.4048,0.8513) -- (3.2529,1.0207) -- (3.0852,1.2487) -- (2.9176,1.5196) -- (2.7500,1.8333) -- (2.8757,2.0968) -- (3.0014,2.3845) -- (3.1271,2.6962) -- (3.2529,3.0321) -- (3.5043,3.7764) -- (3.7461,4.5833) -- cycle;
\fill[blue!3] (1.7539,4.5833) -- (3.7461,4.5833) -- (3.5043,3.7764) -- (3.2529,3.0321) -- (3.1271,2.6962) -- (3.0014,2.3845) -- (2.8757,2.0968) -- (2.7500,1.8333) -- (2.6243,2.0968) -- (2.4986,2.3845) -- (2.3729,2.6962) -- (2.2471,3.0321) -- (1.9957,3.7764) -- (1.7539,4.5833) -- cycle;
\draw[black,line width=.65pt] (0.4443,4.5833) -- (0.5710,4.1497) -- (0.7176,3.6781) -- (0.8643,3.2392) -- (1.0110,2.8333) -- (1.1576,2.4601) -- (1.3043,2.1199) -- (1.4510,1.8125) -- (1.5976,1.5379) -- (1.7233,1.3287) -- (1.8700,1.1152) -- (2.0167,0.9345) -- (2.0795,0.8671) -- (2.1633,0.7866) -- (2.2471,0.7169) -- (2.3100,0.6716) -- (2.3938,0.6207) -- (2.4567,0.5895) -- (2.5405,0.5573) -- (2.6033,0.5402) -- (2.6871,0.5268) -- (2.7500,0.5238) -- (2.8129,0.5268) -- (2.8967,0.5402) -- (2.9595,0.5573) -- (3.0224,0.5805) -- (3.1062,0.6207) -- (3.1690,0.6579) -- (3.2529,0.7169) -- (3.3157,0.7682) -- (3.3995,0.8460) -- (3.4624,0.9113) -- (3.6090,1.0873) -- (3.7557,1.2962) -- (3.9024,1.5379) -- (4.0490,1.8125) -- (4.1957,2.1199) -- (4.3424,2.4601) -- (4.4890,2.8333) -- (4.6357,3.2392) -- (4.7824,3.6781) -- (4.9290,4.1497) -- (5.0557,4.5833);
\draw[blue!65!black,line width=.65pt] (1.7539,4.5833) -- (1.8700,4.1847) -- (1.9957,3.7764) -- (2.1214,3.3922) -- (2.2262,3.0905) -- (2.3519,2.7505) -- (2.4567,2.4857) -- (2.5824,2.1900) -- (2.6871,1.9621) -- (2.8129,1.7106) -- (2.9176,1.5196) -- (3.0224,1.3452) -- (3.1481,1.1582) -- (3.2529,1.0207) -- (3.3786,0.8779) -- (3.5043,0.7592) -- (3.6300,0.6647) -- (3.7348,0.6044) -- (3.7976,0.5762) -- (3.8605,0.5541) -- (3.9233,0.5380) -- (3.9862,0.5279) -- (4.0910,0.5246) -- (4.1957,0.5380) -- (4.2586,0.5541) -- (4.3214,0.5762) -- (4.3843,0.6044) -- (4.4471,0.6385) -- (4.5729,0.7250) -- (4.6986,0.8357) -- (4.8033,0.9463) -- (4.9081,1.0737) -- (5.0338,1.2487) -- (5.1386,1.4129) -- (5.2643,1.6322) -- (5.3900,1.8756) -- (5.5000,2.1084);
\draw[red!65!black,line width=.65pt] (0.0000,2.1084) -- (0.1100,1.8756) -- (0.2357,1.6322) -- (0.3614,1.4129) -- (0.4662,1.2487) -- (0.5710,1.1012) -- (0.6967,0.9463) -- (0.8014,0.8357) -- (0.9271,0.7250) -- (1.0529,0.6385) -- (1.1157,0.6044) -- (1.1786,0.5762) -- (1.2833,0.5427) -- (1.3462,0.5306) -- (1.4090,0.5246) -- (1.4719,0.5246) -- (1.5348,0.5306) -- (1.6395,0.5541) -- (1.7443,0.5943) -- (1.8071,0.6265) -- (1.8700,0.6647) -- (1.9748,0.7418) -- (2.1005,0.8564) -- (2.2052,0.9704) -- (2.3310,1.1293) -- (2.4567,1.3124) -- (2.5824,1.5196) -- (2.6871,1.7106) -- (2.8129,1.9621) -- (2.9176,2.1900) -- (3.0433,2.4857) -- (3.1481,2.7505) -- (3.2738,3.0905) -- (3.3786,3.3922) -- (3.5043,3.7764) -- (3.6300,4.1847) -- (3.7461,4.5833);
\draw[gray!50,thin] (0.0000,0.0000) rectangle (5.5000,4.5833);
\node[font=\small,fill=white,fill opacity=.85,text opacity=1,inner sep=1pt] at (0.4780,0.4795) {$S_{1}$};
\node[font=\small,fill=white,fill opacity=.85,text opacity=1,inner sep=1pt] at (2.7526,0.9612) {$S_{2}$};
\node[font=\small,fill=white,fill opacity=.85,text opacity=1,inner sep=1pt] at (0.5996,2.2404) {$S_{3}$};
\node[font=\small,fill=white,fill opacity=.85,text opacity=1,inner sep=1pt] at (4.9004,2.2404) {$S_{4}$};
\node[font=\small,fill=white,fill opacity=.85,text opacity=1,inner sep=1pt] at (1.2892,3.9437) {$S_{5}$};
\node[font=\small,fill=white,fill opacity=.85,text opacity=1,inner sep=1pt] at (4.2108,3.9437) {$S_{6}$};
\node[font=\small,fill=white,fill opacity=.85,text opacity=1,inner sep=1pt] at (2.7473,3.8448) {$S_{7}$};
\end{tikzpicture}
\par\smallskip (E3)\end{minipage}\hfill
\begin{minipage}{.47\textwidth}\centering
\begin{tikzpicture}[x=1cm,y=1cm,line cap=round,line join=round]
\fill[blue!3] (2.4317,0.0000) -- (1.6583,0.0000) -- (1.6295,0.3512) -- (1.6039,0.6044) -- (1.5696,0.8789) -- (1.5354,1.1006) -- (1.5012,1.2833) -- (1.4584,1.4710) -- (1.4156,1.6248) -- (1.3595,1.7892) -- (1.4794,1.7391) -- (1.5906,1.6833) -- (1.6933,1.6214) -- (1.7874,1.5531) -- (1.8729,1.4783) -- (1.9499,1.3972) -- (2.0184,1.3109) -- (2.0868,1.2067) -- (2.1467,1.0962) -- (2.1980,0.9823) -- (2.2494,0.8452) -- (2.2921,0.7074) -- (2.3349,0.5411) -- (2.3691,0.3813) -- (2.4033,0.1899) -- (2.4317,0.0000) -- cycle;
\fill[blue!6] (1.6583,0.0000) -- (0.0000,0.0000) -- (0.0000,2.0518) -- (0.4271,2.0024) -- (0.7950,1.9426) -- (0.9575,1.9083) -- (1.1030,1.8719) -- (1.2399,1.8313) -- (1.3595,1.7892) -- (1.4156,1.6248) -- (1.4584,1.4710) -- (1.5012,1.2833) -- (1.5354,1.1006) -- (1.5696,0.8789) -- (1.6039,0.6044) -- (1.6295,0.3512) -- (1.6583,0.0000) -- cycle;
\fill[blue!9] (0.0000,2.0518) -- (0.0000,2.7891) -- (0.1409,2.7576) -- (0.2778,2.7218) -- (0.3976,2.6850) -- (0.5088,2.6452) -- (0.6115,2.6023) -- (0.7141,2.5520) -- (0.7997,2.5027) -- (0.8852,2.4451) -- (0.9622,2.3842) -- (1.0313,2.3203) -- (1.0991,2.2466) -- (1.1590,2.1703) -- (1.2189,2.0803) -- (1.2702,1.9893) -- (1.3130,1.9012) -- (1.3595,1.7892) -- (1.2399,1.8313) -- (1.1030,1.8719) -- (0.9575,1.9083) -- (0.7950,1.9426) -- (0.4271,2.0024) -- (0.0000,2.0518) -- cycle;
\fill[blue!3] (5.0703,2.6386) -- (5.0703,0.0000) -- (2.4317,0.0000) -- (2.3948,0.2412) -- (2.3520,0.4647) -- (2.3007,0.6766) -- (2.2494,0.8452) -- (2.1895,1.0028) -- (2.1210,1.1461) -- (2.0440,1.2742) -- (1.9499,1.3972) -- (1.8906,1.4609) -- (1.8302,1.5174) -- (1.7617,1.5730) -- (1.6933,1.6214) -- (1.6163,1.6689) -- (1.5393,1.7103) -- (1.3595,1.7892) -- (1.3130,1.9012) -- (1.2702,1.9893) -- (1.2189,2.0803) -- (1.1590,2.1703) -- (1.0991,2.2466) -- (1.0307,2.3209) -- (0.9537,2.3915) -- (0.8767,2.4513) -- (0.7911,2.5080) -- (0.7056,2.5565) -- (0.6115,2.6023) -- (0.5088,2.6452) -- (0.3976,2.6850) -- (0.2778,2.7218) -- (0.1409,2.7576) -- (0.0000,2.7891) -- (0.0000,5.5000) -- (2.2090,5.5000) -- (2.2544,5.2099) -- (2.2972,4.9963) -- (2.3571,4.7631) -- (2.4169,4.5829) -- (2.4854,4.4214) -- (2.5624,4.2791) -- (2.6479,4.1549) -- (2.7500,4.0391) -- (2.8618,3.9401) -- (2.9816,3.8566) -- (3.1185,3.7812) -- (3.2812,3.7108) -- (3.3276,3.5988) -- (3.3875,3.4787) -- (3.4474,3.3793) -- (3.5159,3.2846) -- (3.5928,3.1965) -- (3.6698,3.1232) -- (3.7554,3.0549) -- (3.8581,2.9868) -- (3.9607,2.9303) -- (4.0805,2.8754) -- (4.2088,2.8266) -- (4.3543,2.7807) -- (4.5082,2.7403) -- (4.6793,2.7031) -- (4.8676,2.6691) -- (5.0703,2.6386) -- cycle;
\fill[blue!6] (5.0703,3.4120) -- (5.0703,2.6386) -- (4.8676,2.6691) -- (4.6793,2.7031) -- (4.5082,2.7403) -- (4.3543,2.7807) -- (4.2088,2.8266) -- (4.0805,2.8754) -- (3.9607,2.9303) -- (3.8581,2.9868) -- (3.7554,3.0549) -- (3.6698,3.1232) -- (3.5928,3.1965) -- (3.5159,3.2846) -- (3.4474,3.3793) -- (3.3875,3.4787) -- (3.3276,3.5988) -- (3.2812,3.7108) -- (3.4350,3.6579) -- (3.5975,3.6124) -- (3.7857,3.5694) -- (3.9911,3.5313) -- (4.2221,3.4965) -- (4.4787,3.4650) -- (4.7610,3.4370) -- (5.0703,3.4120) -- cycle;
\fill[blue!9] (2.9823,5.5000) -- (5.0703,5.5000) -- (5.0703,3.4120) -- (4.7610,3.4370) -- (4.4787,3.4650) -- (4.2221,3.4965) -- (3.9911,3.5313) -- (3.7857,3.5694) -- (3.5975,3.6124) -- (3.4350,3.6579) -- (3.2812,3.7108) -- (3.2250,3.8752) -- (3.1822,4.0290) -- (3.1394,4.2167) -- (3.1052,4.3994) -- (3.0710,4.6211) -- (3.0368,4.8956) -- (3.0111,5.1488) -- (2.9823,5.5000) -- cycle;
\fill[blue!3] (2.2090,5.5000) -- (2.9823,5.5000) -- (3.0111,5.1488) -- (3.0368,4.8956) -- (3.0710,4.6211) -- (3.1052,4.3994) -- (3.1394,4.2167) -- (3.1822,4.0290) -- (3.2250,3.8752) -- (3.2812,3.7108) -- (3.1612,3.7609) -- (3.0500,3.8167) -- (2.9474,3.8786) -- (2.8532,3.9469) -- (2.7677,4.0217) -- (2.6907,4.1028) -- (2.6223,4.1891) -- (2.5538,4.2933) -- (2.4939,4.4038) -- (2.4426,4.5177) -- (2.3913,4.6548) -- (2.3485,4.7926) -- (2.3057,4.9589) -- (2.2715,5.1187) -- (2.2373,5.3101) -- (2.2090,5.5000) -- cycle;
\draw[black,line width=.65pt] (0.0000,2.0518) -- (0.3501,2.0126) -- (0.6496,1.9687) -- (0.9148,1.9179) -- (1.0345,1.8898) -- (1.1457,1.8600) -- (1.2484,1.8285) -- (1.3425,1.7956) -- (1.4366,1.7580) -- (1.5222,1.7188) -- (1.5992,1.6786) -- (1.6762,1.6326) -- (1.7446,1.5858) -- (1.8130,1.5321) -- (1.8729,1.4783) -- (1.9328,1.4166) -- (1.9841,1.3560) -- (2.0355,1.2867) -- (2.0868,1.2067) -- (2.1296,1.1299) -- (2.1724,1.0418) -- (2.2066,0.9613) -- (2.2408,0.8700) -- (2.2750,0.7655) -- (2.3092,0.6447) -- (2.3349,0.5411) -- (2.3862,0.2901) -- (2.4317,0.0000);
\draw[black,line width=.65pt] (2.9823,5.5000) -- (3.0197,5.0591) -- (3.0624,4.6841) -- (3.0881,4.5047) -- (3.1138,4.3505) -- (3.1651,4.0995) -- (3.1908,3.9959) -- (3.2250,3.8752) -- (3.2592,3.7707) -- (3.3020,3.6583) -- (3.3448,3.5620) -- (3.3875,3.4787) -- (3.4303,3.4059) -- (3.4816,3.3297) -- (3.5330,3.2636) -- (3.5928,3.1965) -- (3.6613,3.1307) -- (3.7297,3.0741) -- (3.8067,3.0192) -- (3.8837,2.9717) -- (3.9693,2.9260) -- (4.0548,2.8863) -- (4.1575,2.8450) -- (4.2601,2.8094) -- (4.3714,2.7758) -- (4.4911,2.7445) -- (4.6195,2.7154) -- (4.7563,2.6884) -- (4.9103,2.6622) -- (5.0703,2.6386);
\draw[blue!65!black,line width=.65pt] (0.0000,2.7891) -- (0.1238,2.7617) -- (0.2350,2.7336) -- (0.3377,2.7041) -- (0.4404,2.6705) -- (0.5345,2.6351) -- (0.6200,2.5984) -- (0.6970,2.5610) -- (0.7740,2.5183) -- (0.8425,2.4751) -- (0.9109,2.4259) -- (0.9708,2.3768) -- (1.0307,2.3209) -- (1.0820,2.2664) -- (1.1333,2.2045) -- (1.1847,2.1336) -- (1.2274,2.0661) -- (1.2702,1.9893) -- (1.3130,1.9012) -- (1.3472,1.8207) -- (1.3814,1.7293) -- (1.4156,1.6248) -- (1.4413,1.5360) -- (1.4755,1.4005) -- (1.5012,1.2833) -- (1.5440,1.0493) -- (1.5867,0.7494) -- (1.6210,0.4409) -- (1.6583,0.0000);
\draw[blue!65!black,line width=.65pt] (2.2090,5.5000) -- (2.2544,5.2099) -- (2.2801,5.0761) -- (2.3143,4.9229) -- (2.3485,4.7926) -- (2.3827,4.6805) -- (2.4169,4.5829) -- (2.4597,4.4774) -- (2.5025,4.3867) -- (2.5453,4.3078) -- (2.5966,4.2258) -- (2.6479,4.1549) -- (2.6993,4.0930) -- (2.7591,4.0300) -- (2.8276,3.9679) -- (2.8960,3.9142) -- (2.9730,3.8620) -- (3.0500,3.8167) -- (3.1356,3.7729) -- (3.2297,3.7312) -- (3.3323,3.6919) -- (3.4350,3.6579) -- (3.5462,3.6258) -- (3.6745,3.5937) -- (3.8029,3.5659) -- (3.9483,3.5386) -- (4.1023,3.5136) -- (4.2734,3.4896) -- (4.6413,3.4482) -- (5.0703,3.4120);
\draw[gray!50,thin] (0.0000,0.0000) rectangle (5.0703,5.5000);
\node[font=\small,fill=white,fill opacity=.85,text opacity=1,inner sep=1pt] at (1.9969,0.3669) {$S_{1}$};
\node[font=\small,fill=white,fill opacity=.85,text opacity=1,inner sep=1pt] at (0.7871,0.7874) {$S_{2}$};
\node[font=\small,fill=white,fill opacity=.85,text opacity=1,inner sep=1pt] at (0.3385,2.3515) {$S_{3}$};
\node[font=\small,fill=white,fill opacity=.85,text opacity=1,inner sep=1pt] at (3.5919,1.4784) {$S_{4}$};
\node[font=\small,fill=white,fill opacity=.85,text opacity=1,inner sep=1pt] at (4.7034,3.0734) {$S_{5}$};
\node[font=\small,fill=white,fill opacity=.85,text opacity=1,inner sep=1pt] at (4.0855,4.5151) {$S_{6}$};
\node[font=\small,fill=white,fill opacity=.85,text opacity=1,inner sep=1pt] at (2.6437,5.1331) {$S_{7}$};
\end{tikzpicture}
\par\smallskip (E4)\end{minipage}
\caption{Three parabolas (E3) and two hyperbolas (E4).}\label{fig:second}
\end{figure}
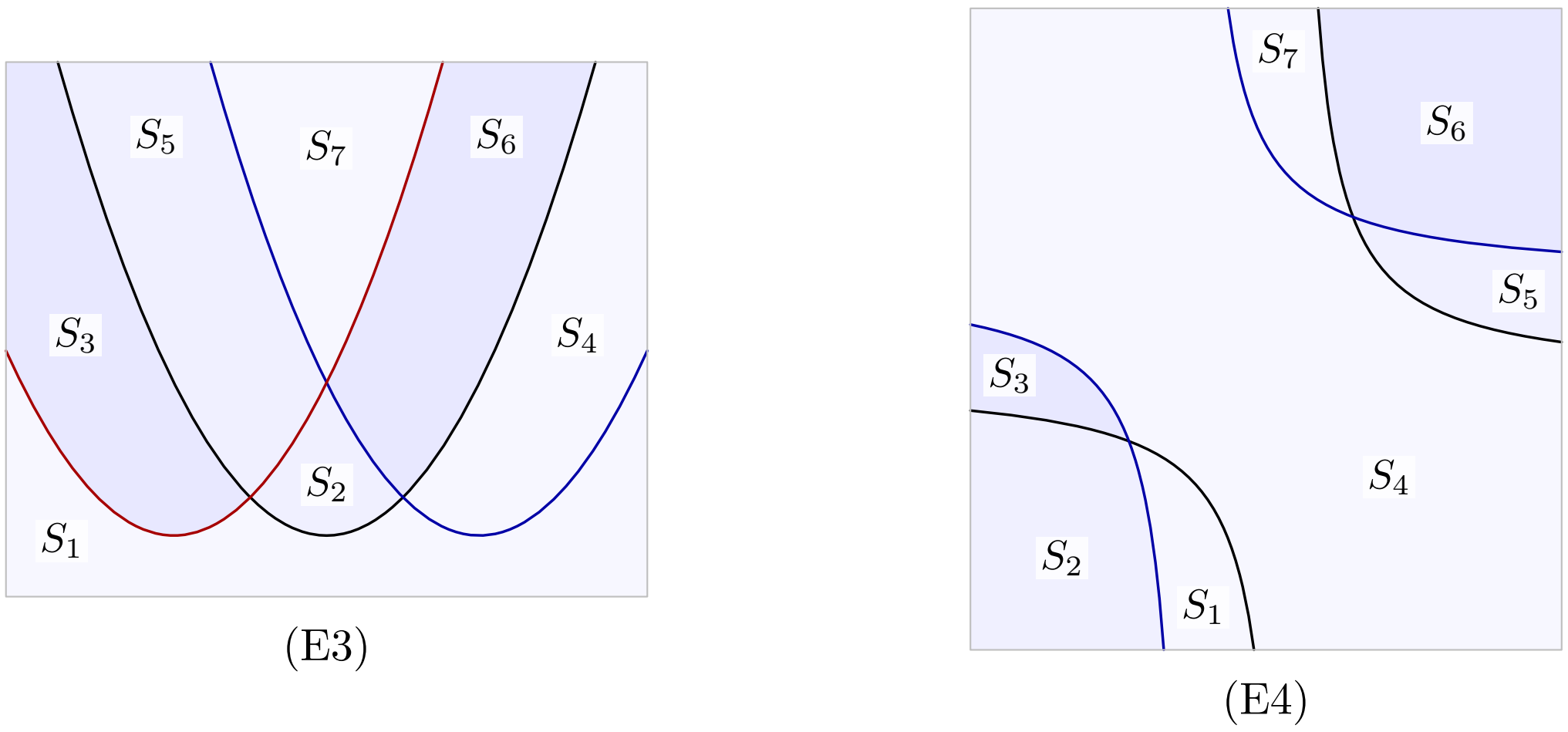
\clearpage
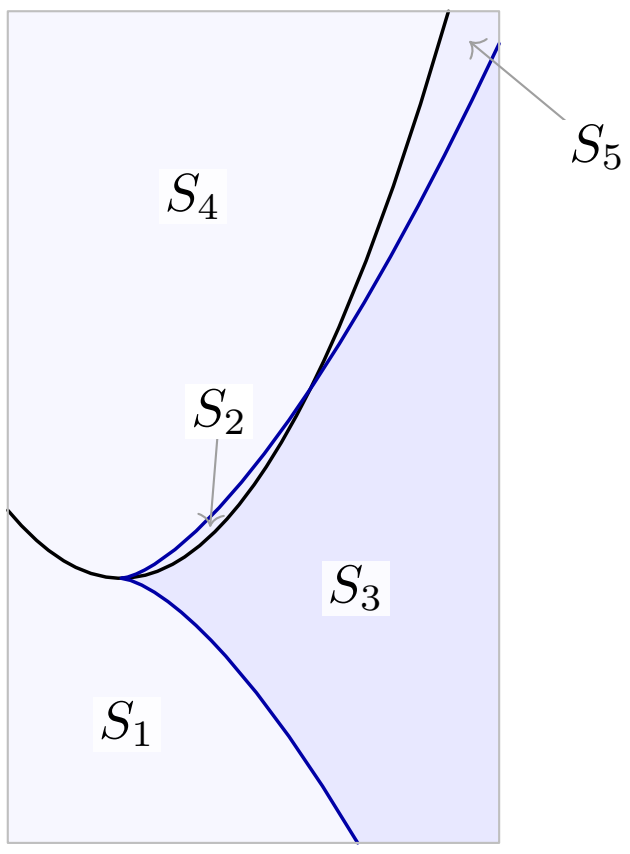
\begin{figure}[htbp]\centering
\begin{minipage}{.47\textwidth}\centering
\begin{tikzpicture}[x=1cm,y=1cm,line cap=round,line join=round]
\fill[blue!3] (4.1532,0.0000) -- (0.0000,0.0000) -- (0.0000,4.6667) -- (4.1532,4.6667) -- (3.7110,3.9246) -- (3.3059,3.3149) -- (3.0915,3.0218) -- (2.8868,2.7627) -- (2.6919,2.5360) -- (2.5000,2.3333) -- (2.6919,2.1307) -- (2.8868,1.9039) -- (3.0915,1.6448) -- (3.3059,1.3518) -- (3.7110,0.7421) -- (4.1532,0.0000) -- cycle;
\fill[blue!6] (4.1532,4.6667) -- (5.5000,4.6667) -- (5.5000,0.0000) -- (4.1532,0.0000) -- (3.7110,0.7421) -- (3.3059,1.3518) -- (3.0915,1.6448) -- (2.8868,1.9039) -- (2.6919,2.1307) -- (2.5000,2.3333) -- (2.6919,2.5360) -- (2.8868,2.7627) -- (3.0915,3.0218) -- (3.3059,3.3149) -- (3.7110,3.9246) -- (4.1532,4.6667) -- cycle;
\fill[blue!9] (2.5000,2.3333) -- (2.3312,2.1733) -- (2.1833,2.0483) -- (2.0432,1.9441) -- (1.8947,1.8503) -- (1.7559,1.7797) -- (1.6268,1.7308) -- (1.5075,1.7021) -- (1.3863,1.6918) -- (1.2668,1.7044) -- (1.1619,1.7392) -- (1.0715,1.7934) -- (0.9894,1.8711) -- (0.9246,1.9647) -- (0.8739,2.0797) -- (0.8435,2.2041) -- (0.8333,2.3333) -- (0.8435,2.4625) -- (0.8739,2.5870) -- (0.9246,2.7020) -- (0.9894,2.7956) -- (1.0715,2.8733) -- (1.1619,2.9275) -- (1.2668,2.9622) -- (1.3863,2.9748) -- (1.5075,2.9646) -- (1.6268,2.9358) -- (1.7559,2.8869) -- (1.8947,2.8164) -- (2.0432,2.7225) -- (2.1833,2.6183) -- (2.3312,2.4934) -- (2.5000,2.3333) -- cycle;
\draw[black,line width=.65pt] (4.1532,0.0000) -- (3.9251,0.3924) -- (3.7110,0.7421) -- (3.5046,1.0615) -- (3.3059,1.3518) -- (3.0915,1.6448) -- (2.8868,1.9039) -- (2.6919,2.1307) -- (2.5000,2.3333);
\draw[black,line width=.65pt] (2.5000,2.3333) -- (2.3312,2.4934) -- (2.1833,2.6183) -- (2.0432,2.7225) -- (1.8947,2.8164) -- (1.7559,2.8869) -- (1.6268,2.9358) -- (1.5075,2.9646) -- (1.3863,2.9748) -- (1.2668,2.9622) -- (1.1619,2.9275) -- (1.0715,2.8733) -- (0.9894,2.7956) -- (0.9246,2.7020) -- (0.8739,2.5870) -- (0.8435,2.4625) -- (0.8333,2.3333) -- (0.8435,2.2041) -- (0.8739,2.0797) -- (0.9246,1.9647) -- (0.9894,1.8711) -- (1.0715,1.7934) -- (1.1619,1.7392) -- (1.2668,1.7044) -- (1.3863,1.6918) -- (1.5075,1.7021) -- (1.6268,1.7308) -- (1.7559,1.7797) -- (1.8947,1.8503) -- (2.0432,1.9441) -- (2.1833,2.0483) -- (2.3312,2.1733) -- (2.5000,2.3333);
\draw[black,line width=.65pt] (2.5000,2.3333) -- (2.6919,2.5360) -- (2.8868,2.7627) -- (3.0915,3.0218) -- (3.3059,3.3149) -- (3.5046,3.6051) -- (3.7110,3.9246) -- (3.9251,4.2743) -- (4.1532,4.6667);
\draw[gray!50,thin] (0.0000,0.0000) rectangle (5.5000,4.6667);
\node[font=\small,fill=white,fill opacity=.85,text opacity=1,inner sep=1pt] at (2.3569,0.9821) {$S_{1}$};
\node[font=\small,fill=white,fill opacity=.85,text opacity=1,inner sep=1pt] at (4.2028,2.3265) {$S_{2}$};
\node[font=\small,fill=white,fill opacity=.85,text opacity=1,inner sep=1pt] at (1.4551,2.3366) {$S_{3}$};
\end{tikzpicture}
\par\smallskip (E5)\end{minipage}\hfill
\begin{minipage}{.47\textwidth}\centering
\begin{tikzpicture}[x=1cm,y=1cm,line cap=round,line join=round]
\fill[blue!3] (2.3143,0.0000) -- (0.0000,0.0000) -- (0.0000,2.2001) -- (0.0900,2.0985) -- (0.1900,2.0009) -- (0.2700,1.9343) -- (0.3700,1.8655) -- (0.4500,1.8220) -- (0.5500,1.7820) -- (0.6500,1.7580) -- (0.7781,1.7458) -- (0.8861,1.7051) -- (0.9705,1.6574) -- (1.0600,1.5956) -- (1.1561,1.5185) -- (1.2556,1.4284) -- (1.3766,1.3064) -- (1.4988,1.1704) -- (1.6320,1.0091) -- (1.7625,0.8388) -- (2.0352,0.4468) -- (2.3143,0.0000) -- cycle;
\fill[blue!6] (2.0000,3.0000) -- (1.8500,2.7180) -- (1.6900,2.4569) -- (1.5300,2.2367) -- (1.3700,2.0575) -- (1.2300,1.9343) -- (1.1500,1.8780) -- (1.0700,1.8319) -- (0.9900,1.7961) -- (0.9100,1.7705) -- (0.8300,1.7551) -- (0.7900,1.7513) -- (0.7700,1.7503) -- (0.7503,1.7500) -- (0.7990,1.7597) -- (0.8580,1.7818) -- (0.9964,1.8594) -- (1.1476,1.9743) -- (1.3044,2.1193) -- (1.4758,2.3031) -- (1.6446,2.5069) -- (1.8172,2.7361) -- (2.0000,3.0000) -- cycle;
\fill[blue!9] (3.2500,5.2856) -- (3.2500,0.0000) -- (2.3143,0.0000) -- (2.0813,0.3761) -- (1.8592,0.7051) -- (1.6446,0.9931) -- (1.4419,1.2352) -- (1.2461,1.4374) -- (1.0600,1.5956) -- (0.9705,1.6574) -- (0.8861,1.7051) -- (0.8116,1.7363) -- (0.7500,1.7500) -- (0.8500,1.7580) -- (0.9700,1.7887) -- (1.0700,1.8319) -- (1.1900,1.9049) -- (1.3100,2.0009) -- (1.4100,2.0985) -- (1.5300,2.2367) -- (1.6300,2.3695) -- (1.7900,2.6153) -- (1.8900,2.7897) -- (2.0000,3.0000) -- (2.2735,3.4320) -- (2.5500,3.9100) -- (2.8886,4.5473) -- (3.2500,5.2856) -- cycle;
\fill[blue!3] (0.0000,2.2001) -- (0.0000,5.5000) -- (2.9150,5.5000) -- (2.7300,4.8863) -- (2.5300,4.2847) -- (2.3300,3.7471) -- (2.1300,3.2735) -- (2.0000,3.0000) -- (1.9603,2.9410) -- (1.7092,2.5903) -- (1.4644,2.2901) -- (1.2556,2.0716) -- (1.0675,1.9100) -- (0.9768,1.8466) -- (0.8911,1.7974) -- (0.8150,1.7648) -- (0.7522,1.7501) -- (0.6700,1.7551) -- (0.5700,1.7759) -- (0.4700,1.8127) -- (0.3700,1.8655) -- (0.2700,1.9343) -- (0.1900,2.0009) -- (0.0900,2.0985) -- (0.0000,2.2001) -- cycle;
\fill[blue!6] (2.9150,5.5000) -- (3.2500,5.5000) -- (3.2500,5.2856) -- (2.9280,4.6250) -- (2.6044,4.0087) -- (2.2901,3.4595) -- (2.0000,3.0000) -- (2.1100,3.2297) -- (2.2300,3.5023) -- (2.4500,4.0620) -- (2.6900,4.7609) -- (2.9150,5.5000) -- cycle;
\draw[black,line width=.65pt] (0.0000,2.2001) -- (0.0900,2.0985) -- (0.1900,2.0009) -- (0.2700,1.9343) -- (0.3700,1.8655) -- (0.4500,1.8220) -- (0.5500,1.7820) -- (0.6500,1.7580) -- (0.7300,1.7503) -- (0.8100,1.7529) -- (0.9100,1.7705) -- (0.9900,1.7961) -- (1.0900,1.8425) -- (1.1700,1.8911) -- (1.2700,1.9663) -- (1.3700,2.0575) -- (1.4500,2.1420) -- (1.5300,2.2367) -- (1.6300,2.3695) -- (1.7100,2.4873) -- (1.8100,2.6489) -- (1.9100,2.8265) -- (2.0000,3.0000) -- (2.0900,3.1865) -- (2.1900,3.4089) -- (2.3700,3.8495) -- (2.5500,4.3420) -- (2.7300,4.8863) -- (2.9150,5.5000);
\draw[blue!65!black,line width=.65pt] (2.3143,0.0000) -- (2.0813,0.3761) -- (1.8592,0.7051) -- (1.6446,0.9931) -- (1.4419,1.2352) -- (1.3451,1.3394) -- (1.2461,1.4374) -- (1.1476,1.5257) -- (1.0600,1.5956) -- (0.9705,1.6574) -- (0.8861,1.7051) -- (0.8116,1.7363) -- (0.7500,1.7500) -- (0.7932,1.7580) -- (0.8494,1.7780) -- (0.9066,1.8055) -- (0.9705,1.8426) -- (1.1064,1.9403) -- (1.2461,2.0626) -- (1.3980,2.2166) -- (1.5460,2.3852) -- (1.6961,2.5731) -- (1.8592,2.7949) -- (2.0201,3.0302) -- (2.1918,3.2985) -- (2.3574,3.5728) -- (2.5320,3.8778) -- (2.7156,4.2149) -- (2.8886,4.5473) -- (3.0688,4.9082) -- (3.2500,5.2856);
\draw[gray!50,thin] (0.0000,0.0000) rectangle (3.2500,5.5000);
\node[font=\small,fill=white,fill opacity=.85,text opacity=1,inner sep=1pt] at (0.7888,0.7859) {$S_{1}$};
\draw[gray!75,thin,->] (1.4000,2.8500) -- (1.3385,2.0894);
\node[font=\small,fill=white,inner sep=1pt] at (1.4000,2.8500) {$S_{2}$};
\node[font=\small,fill=white,fill opacity=.85,text opacity=1,inner sep=1pt] at (2.3020,1.6853) {$S_{3}$};
\node[font=\small,fill=white,fill opacity=.85,text opacity=1,inner sep=1pt] at (1.2250,4.2750) {$S_{4}$};
\draw[gray!75,thin,->] (3.9000,4.6000) -- (3.0522,5.3022);
\node[font=\small,fill=white,inner sep=1pt] at (3.9000,4.6000) {$S_{5}$};
\end{tikzpicture}
\par\smallskip (E6)\end{minipage}
\caption{A nodal cubic (E5) and a parabola with a cuspidal cubic (E6).
In (E6), arrows point to the narrow regions $S_2$ and $S_5$.}\label{fig:third}
\end{figure}
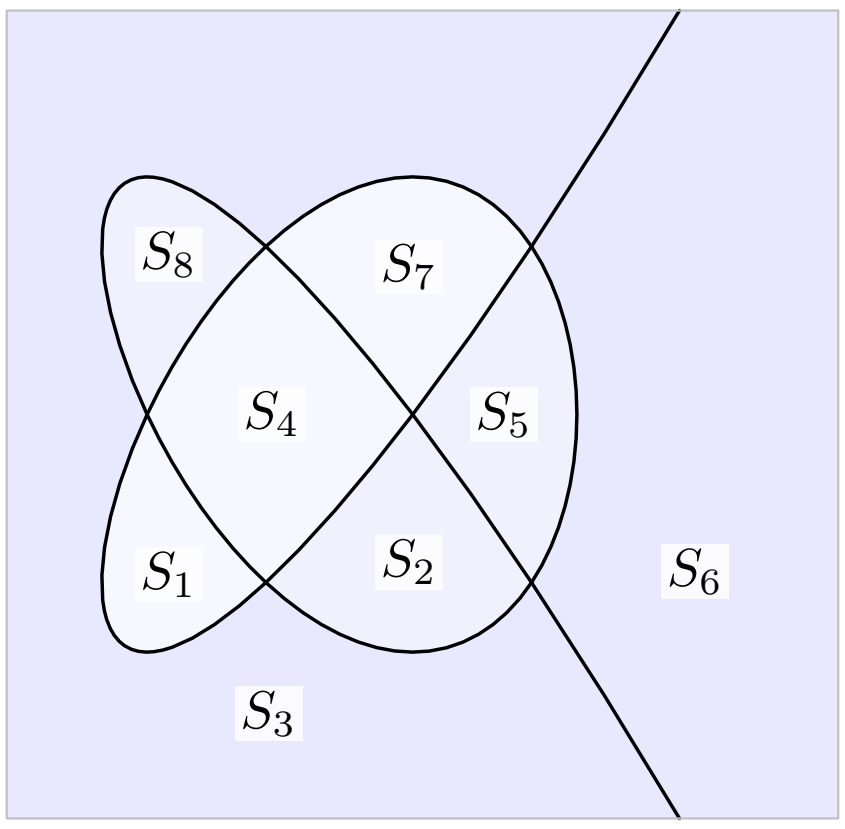
\begin{figure}[htbp]\centering
\begin{minipage}{.47\textwidth}\centering
\begin{tikzpicture}[x=1cm,y=1cm,line cap=round,line join=round]
\fill[blue!3] (1.5603,2.9857) -- (1.8573,2.5161) -- (2.1386,2.1020) -- (2.4085,1.7364) -- (2.6714,1.4143) -- (2.5032,1.2284) -- (2.3477,1.0724) -- (2.1959,0.9370) -- (2.0552,0.8285) -- (1.9197,0.7419) -- (1.7906,0.6789) -- (1.6692,0.6408) -- (1.5607,0.6286) -- (1.4575,0.6418) -- (1.3685,0.6799) -- (1.2881,0.7455) -- (1.2228,0.8337) -- (1.1690,0.9496) -- (1.1311,1.0844) -- (1.1079,1.2399) -- (1.1000,1.4139) -- (1.1080,1.6039) -- (1.1314,1.8011) -- (1.1706,2.0079) -- (1.2240,2.2151) -- (1.2906,2.4192) -- (1.3694,2.6172) -- (1.4593,2.8063) -- (1.5603,2.9857) -- cycle;
\fill[blue!6] (3.7826,2.9857) -- (3.8836,2.8063) -- (3.9735,2.6172) -- (4.0522,2.4192) -- (4.1189,2.2151) -- (4.1723,2.0079) -- (4.2115,1.8011) -- (4.2349,1.6039) -- (4.2429,1.4139) -- (4.2350,1.2399) -- (4.2117,1.0844) -- (4.1739,0.9496) -- (4.1201,0.8337) -- (4.0548,0.7455) -- (3.9743,0.6799) -- (3.8853,0.6418) -- (3.7822,0.6286) -- (3.6736,0.6408) -- (3.5522,0.6789) -- (3.4231,0.7419) -- (3.2877,0.8285) -- (3.1470,0.9370) -- (2.9952,1.0724) -- (2.8396,1.2284) -- (2.6714,1.4143) -- (2.9343,1.7364) -- (3.2042,2.1020) -- (3.4856,2.5161) -- (3.7826,2.9857) -- cycle;
\fill[blue!9] (2.6714,1.4143) -- (2.4085,1.7364) -- (2.1386,2.1020) -- (1.8573,2.5161) -- (1.5603,2.9857) -- (1.6778,3.1611) -- (1.8057,3.3194) -- (1.9415,3.4566) -- (2.0800,3.5680) -- (2.2235,3.6562) -- (2.3707,3.7199) -- (2.5205,3.7585) -- (2.6714,3.7714) -- (2.8224,3.7585) -- (2.9721,3.7199) -- (3.1193,3.6562) -- (3.2628,3.5680) -- (3.4014,3.4566) -- (3.5372,3.3194) -- (3.6651,3.1611) -- (3.7826,2.9857) -- (3.4856,2.5161) -- (3.2042,2.1020) -- (2.9343,1.7364) -- (2.6714,1.4143) -- cycle;
\fill[blue!3] (5.1521,5.5000) -- (5.3429,5.5000) -- (5.3429,0.0000) -- (0.0000,0.0000) -- (0.0000,5.5000) -- (0.1908,5.5000) -- (0.5444,4.8023) -- (0.8983,4.1359) -- (1.2313,3.5406) -- (1.5603,2.9857) -- (1.4593,2.8063) -- (1.3694,2.6172) -- (1.2906,2.4192) -- (1.2240,2.2151) -- (1.1706,2.0079) -- (1.1314,1.8011) -- (1.1080,1.6039) -- (1.1000,1.4139) -- (1.1079,1.2399) -- (1.1311,1.0844) -- (1.1690,0.9496) -- (1.2228,0.8337) -- (1.2881,0.7455) -- (1.3685,0.6799) -- (1.4575,0.6418) -- (1.5607,0.6286) -- (1.6692,0.6408) -- (1.7906,0.6789) -- (1.9197,0.7419) -- (2.0552,0.8285) -- (2.1959,0.9370) -- (2.3477,1.0724) -- (2.5032,1.2284) -- (2.6714,1.4143) -- (2.9076,1.1582) -- (3.1269,0.9539) -- (3.3361,0.7953) -- (3.5263,0.6899) -- (3.6170,0.6558) -- (3.7008,0.6356) -- (3.7781,0.6286) -- (3.8525,0.6345) -- (3.9199,0.6531) -- (3.9834,0.6855) -- (4.0394,0.7300) -- (4.0904,0.7885) -- (4.1374,0.8652) -- (4.1768,0.9578) -- (4.2077,1.0661) -- (4.2287,1.1848) -- (4.2405,1.3174) -- (4.2424,1.4577) -- (4.2175,1.7601) -- (4.1555,2.0791) -- (4.0585,2.4019) -- (3.9337,2.7049) -- (3.7826,2.9857) -- (4.1116,3.5406) -- (4.4446,4.1359) -- (4.7985,4.8023) -- (5.1521,5.5000) -- cycle;
\fill[blue!6] (0.1908,5.5000) -- (5.1521,5.5000) -- (4.7985,4.8023) -- (4.4446,4.1359) -- (4.1116,3.5406) -- (3.7826,2.9857) -- (3.6651,3.1611) -- (3.5372,3.3194) -- (3.4014,3.4566) -- (3.2628,3.5680) -- (3.1193,3.6562) -- (2.9721,3.7199) -- (2.8224,3.7585) -- (2.6714,3.7714) -- (2.5205,3.7585) -- (2.3707,3.7199) -- (2.2235,3.6562) -- (2.0800,3.5680) -- (1.9415,3.4566) -- (1.8057,3.3194) -- (1.6778,3.1611) -- (1.5603,2.9857) -- (1.2313,3.5406) -- (0.8983,4.1359) -- (0.5444,4.8023) -- (0.1908,5.5000) -- cycle;
\draw[black,line width=.65pt] (0.1908,5.5000) -- (0.5444,4.8023) -- (0.8983,4.1359) -- (1.2313,3.5406) -- (1.5603,2.9857);
\draw[black,line width=.65pt] (1.5603,2.9857) -- (1.8573,2.5161) -- (2.1386,2.1020) -- (2.4085,1.7364) -- (2.6714,1.4143);
\draw[black,line width=.65pt] (2.6714,1.4143) -- (2.9076,1.1582) -- (3.0236,1.0457) -- (3.1336,0.9482) -- (3.2378,0.8649) -- (3.3421,0.7914) -- (3.4400,0.7326) -- (3.5316,0.6876) -- (3.6218,0.6543) -- (3.7053,0.6349) -- (3.7822,0.6286) -- (3.8562,0.6352) -- (3.9266,0.6558) -- (3.9893,0.6894) -- (4.0446,0.7350) -- (4.0949,0.7947) -- (4.1410,0.8724) -- (4.1796,0.9660) -- (4.2098,1.0752) -- (4.2300,1.1947) -- (4.2408,1.3227) -- (4.2423,1.4633) -- (4.2344,1.6096) -- (4.2167,1.7659) -- (4.1897,1.9250) -- (4.1541,2.0850) -- (4.1101,2.2445) -- (4.0585,2.4019) -- (3.9995,2.5558) -- (3.9337,2.7049) -- (3.8616,2.8480) -- (3.7826,2.9857);
\draw[black,line width=.65pt] (3.7826,2.9857) -- (3.6651,3.1611) -- (3.5372,3.3194) -- (3.4014,3.4566) -- (3.2628,3.5680) -- (3.1193,3.6562) -- (2.9721,3.7199) -- (2.8224,3.7585) -- (2.7449,3.7684) -- (2.6714,3.7714) -- (2.5979,3.7684) -- (2.5205,3.7585) -- (2.3707,3.7199) -- (2.2235,3.6562) -- (2.0800,3.5680) -- (1.9415,3.4566) -- (1.8057,3.3194) -- (1.6778,3.1611) -- (1.5603,2.9857);
\draw[black,line width=.65pt] (1.5603,2.9857) -- (1.4812,2.8480) -- (1.4092,2.7049) -- (1.3434,2.5558) -- (1.2844,2.4019) -- (1.2327,2.2445) -- (1.1888,2.0850) -- (1.1531,1.9250) -- (1.1262,1.7659) -- (1.1085,1.6096) -- (1.1006,1.4633) -- (1.1021,1.3227) -- (1.1128,1.1947) -- (1.1331,1.0752) -- (1.1632,0.9660) -- (1.2018,0.8724) -- (1.2480,0.7947) -- (1.2983,0.7350) -- (1.3535,0.6894) -- (1.4163,0.6558) -- (1.4866,0.6352) -- (1.5607,0.6286) -- (1.6376,0.6349) -- (1.7210,0.6543) -- (1.8113,0.6876) -- (1.9029,0.7326) -- (2.0008,0.7914) -- (2.1051,0.8649) -- (2.2092,0.9482) -- (2.3193,1.0457) -- (2.4353,1.1582) -- (2.6714,1.4143);
\draw[black,line width=.65pt] (2.6714,1.4143) -- (2.9343,1.7364) -- (3.2042,2.1020) -- (3.4856,2.5161) -- (3.7826,2.9857);
\draw[black,line width=.65pt] (3.7826,2.9857) -- (4.1116,3.5406) -- (4.4446,4.1359) -- (4.7985,4.8023) -- (5.1521,5.5000);
\draw[gray!50,thin] (0.0000,0.0000) rectangle (5.3429,5.5000);
\node[font=\small,fill=white,fill opacity=.85,text opacity=1,inner sep=1pt] at (1.7729,1.4596) {$S_{1}$};
\node[font=\small,fill=white,fill opacity=.85,text opacity=1,inner sep=1pt] at (3.5699,1.4596) {$S_{2}$};
\node[font=\small,fill=white,fill opacity=.85,text opacity=1,inner sep=1pt] at (2.6737,2.8898) {$S_{3}$};
\node[font=\small,fill=white,fill opacity=.85,text opacity=1,inner sep=1pt] at (4.6185,2.9514) {$S_{4}$};
\node[font=\small,fill=white,fill opacity=.85,text opacity=1,inner sep=1pt] at (1.7945,4.5287) {$S_{5}$};
\end{tikzpicture}
\par\smallskip (E7)\end{minipage}\hfill
\begin{minipage}{.47\textwidth}\centering
\begin{tikzpicture}[x=1cm,y=1cm,line cap=round,line join=round]
\fill[blue!3] (0.9287,2.6714) -- (1.0937,2.3565) -- (1.2861,2.0562) -- (1.4966,1.7868) -- (1.7144,1.5603) -- (1.4794,1.3592) -- (1.2695,1.2156) -- (1.1734,1.1648) -- (1.0861,1.1293) -- (1.0013,1.1069) -- (0.9261,1.1000) -- (0.8602,1.1076) -- (0.7993,1.1305) -- (0.7484,1.1675) -- (0.7045,1.2212) -- (0.6710,1.2880) -- (0.6475,1.3665) -- (0.6336,1.4552) -- (0.6286,1.5595) -- (0.6482,1.8002) -- (0.7063,2.0731) -- (0.8009,2.3677) -- (0.9287,2.6714) -- cycle;
\fill[blue!6] (2.6856,2.6714) -- (3.0681,2.1491) -- (3.4713,1.5603) -- (3.3944,1.4560) -- (3.3076,1.3614) -- (3.2175,1.2832) -- (3.1208,1.2179) -- (3.0179,1.1663) -- (2.9094,1.1291) -- (2.7958,1.1068) -- (2.6832,1.1000) -- (2.5670,1.1075) -- (2.4479,1.1297) -- (2.3265,1.1665) -- (2.2033,1.2180) -- (2.0791,1.2841) -- (1.9605,1.3604) -- (1.8364,1.4541) -- (1.7144,1.5603) -- (1.9304,1.7735) -- (2.1715,2.0373) -- (2.4222,2.3356) -- (2.6856,2.6714) -- cycle;
\fill[blue!9] (4.4497,0.0000) -- (0.0000,0.0000) -- (0.0000,5.3429) -- (4.4497,5.3429) -- (3.9523,4.5284) -- (3.4713,3.7826) -- (3.3944,3.8869) -- (3.3076,3.9815) -- (3.2175,4.0597) -- (3.1208,4.1249) -- (3.0179,4.1765) -- (2.9094,4.2138) -- (2.7958,4.2360) -- (2.6832,4.2429) -- (2.5670,4.2353) -- (2.4479,4.2132) -- (2.3265,4.1763) -- (2.2033,4.1248) -- (2.0791,4.0587) -- (1.9605,3.9825) -- (1.8364,3.8888) -- (1.7144,3.7826) -- (1.4598,3.9987) -- (1.3386,4.0844) -- (1.2284,4.1503) -- (1.1214,4.2006) -- (1.0263,4.2308) -- (0.9426,4.2426) -- (0.8649,4.2364) -- (0.8115,4.2185) -- (0.7655,4.1901) -- (0.7236,4.1485) -- (0.6899,4.0964) -- (0.6638,4.0350) -- (0.6440,3.9594) -- (0.6286,3.7833) -- (0.6482,3.5427) -- (0.7063,3.2698) -- (0.8009,2.9751) -- (0.9287,2.6714) -- (0.8236,2.4273) -- (0.7402,2.1899) -- (0.6780,1.9584) -- (0.6408,1.7467) -- (0.6286,1.5595) -- (0.6419,1.3950) -- (0.6810,1.2644) -- (0.7097,1.2132) -- (0.7451,1.1706) -- (0.8157,1.1224) -- (0.9048,1.1008) -- (1.0075,1.1080) -- (1.1287,1.1452) -- (1.2611,1.2107) -- (1.4024,1.3020) -- (1.5502,1.4159) -- (1.7144,1.5603) -- (1.8541,1.4398) -- (1.9901,1.3401) -- (2.1324,1.2540) -- (2.2738,1.1868) -- (2.4134,1.1387) -- (2.5502,1.1098) -- (2.6832,1.1000) -- (2.8116,1.1090) -- (2.9043,1.1278) -- (2.9986,1.1585) -- (3.0887,1.2000) -- (3.1743,1.2519) -- (3.2552,1.3137) -- (3.3310,1.3848) -- (3.4051,1.4693) -- (3.4713,1.5603) -- (3.9523,0.8145) -- (4.4497,0.0000) -- cycle;
\fill[blue!3] (1.7144,3.7826) -- (1.9304,3.5693) -- (2.1715,3.3056) -- (2.4222,3.0073) -- (2.6856,2.6714) -- (2.4222,2.3356) -- (2.1715,2.0373) -- (1.9304,1.7735) -- (1.7144,1.5603) -- (1.4966,1.7868) -- (1.2861,2.0562) -- (1.0937,2.3565) -- (0.9287,2.6714) -- (1.0937,2.9864) -- (1.2861,3.2866) -- (1.4966,3.5561) -- (1.7144,3.7826) -- cycle;
\fill[blue!6] (3.4713,1.5603) -- (3.0681,2.1491) -- (2.6856,2.6714) -- (3.0681,3.1938) -- (3.4713,3.7826) -- (3.5401,3.6678) -- (3.5992,3.5466) -- (3.6499,3.4176) -- (3.6937,3.2758) -- (3.7266,3.1346) -- (3.7515,2.9827) -- (3.7664,2.8278) -- (3.7714,2.6714) -- (3.7664,2.5151) -- (3.7515,2.3602) -- (3.7266,2.2083) -- (3.6937,2.0671) -- (3.6499,1.9252) -- (3.5992,1.7962) -- (3.5401,1.6750) -- (3.4713,1.5603) -- cycle;
\fill[blue!9] (4.4497,5.3429) -- (5.5000,5.3429) -- (5.5000,0.0000) -- (4.4497,0.0000) -- (3.9523,0.8145) -- (3.4713,1.5603) -- (3.5401,1.6750) -- (3.5992,1.7962) -- (3.6499,1.9252) -- (3.6937,2.0671) -- (3.7266,2.2083) -- (3.7515,2.3602) -- (3.7664,2.5151) -- (3.7714,2.6714) -- (3.7664,2.8278) -- (3.7515,2.9827) -- (3.7266,3.1346) -- (3.6937,3.2758) -- (3.6499,3.4176) -- (3.5992,3.5466) -- (3.5401,3.6678) -- (3.4713,3.7826) -- (3.9523,4.5284) -- (4.4497,5.3429) -- cycle;
\fill[blue!3] (2.6856,2.6714) -- (2.4222,3.0073) -- (2.1715,3.3056) -- (1.9304,3.5693) -- (1.7144,3.7826) -- (1.8364,3.8888) -- (1.9605,3.9825) -- (2.0791,4.0587) -- (2.2033,4.1248) -- (2.3265,4.1763) -- (2.4479,4.2132) -- (2.5670,4.2353) -- (2.6832,4.2429) -- (2.7958,4.2360) -- (2.9094,4.2138) -- (3.0179,4.1765) -- (3.1208,4.1249) -- (3.2175,4.0597) -- (3.3076,3.9815) -- (3.3944,3.8869) -- (3.4713,3.7826) -- (3.0681,3.1938) -- (2.6856,2.6714) -- cycle;
\fill[blue!6] (0.9287,2.6714) -- (0.8009,2.9751) -- (0.7063,3.2698) -- (0.6482,3.5427) -- (0.6286,3.7833) -- (0.6336,3.8877) -- (0.6475,3.9764) -- (0.6710,4.0548) -- (0.7045,4.1217) -- (0.7484,4.1754) -- (0.7993,4.2124) -- (0.8602,4.2353) -- (0.9261,4.2428) -- (1.0013,4.2360) -- (1.0861,4.2136) -- (1.1734,4.1781) -- (1.2695,4.1273) -- (1.4794,3.9836) -- (1.7144,3.7826) -- (1.4966,3.5561) -- (1.2861,3.2866) -- (1.0937,2.9864) -- (0.9287,2.6714) -- cycle;
\draw[black,line width=.65pt] (4.4497,0.0000) -- (3.9523,0.8145) -- (3.4713,1.5603);
\draw[black,line width=.65pt] (3.4713,1.5603) -- (3.0681,2.1491) -- (2.6856,2.6714);
\draw[black,line width=.65pt] (2.6856,2.6714) -- (2.4222,3.0073) -- (2.1715,3.3056) -- (1.9304,3.5693) -- (1.7144,3.7826);
\draw[black,line width=.65pt] (1.7144,3.7826) -- (1.5399,3.9354) -- (1.3747,4.0602) -- (1.2284,4.1503) -- (1.0931,4.2112) -- (1.0327,4.2293) -- (0.9712,4.2404) -- (0.9207,4.2428) -- (0.8697,4.2374) -- (0.8242,4.2240) -- (0.7800,4.2007) -- (0.7451,4.1722) -- (0.7124,4.1335) -- (0.6854,4.0876) -- (0.6621,4.0299) -- (0.6463,3.9708) -- (0.6349,3.9001) -- (0.6295,3.7349) -- (0.6470,3.5504) -- (0.6875,3.3438) -- (0.7482,3.1277) -- (0.8303,2.8985) -- (0.9287,2.6714);
\draw[black,line width=.65pt] (0.9287,2.6714) -- (1.0083,2.5112) -- (1.0937,2.3565) -- (1.1848,2.2069) -- (1.2861,2.0562) -- (1.3868,1.9206) -- (1.4966,1.7868) -- (1.6039,1.6690) -- (1.7144,1.5603);
\draw[black,line width=.65pt] (1.7144,1.5603) -- (1.8364,1.4541) -- (1.9605,1.3604) -- (2.0791,1.2841) -- (2.2033,1.2180) -- (2.3265,1.1665) -- (2.4479,1.1297) -- (2.5670,1.1075) -- (2.6832,1.1000) -- (2.7958,1.1068) -- (2.9094,1.1291) -- (3.0179,1.1663) -- (3.1208,1.2179) -- (3.2175,1.2832) -- (3.3076,1.3614) -- (3.3944,1.4560) -- (3.4713,1.5603);
\draw[black,line width=.65pt] (3.4713,1.5603) -- (3.5401,1.6750) -- (3.5992,1.7962) -- (3.6499,1.9252) -- (3.6937,2.0671) -- (3.7266,2.2083) -- (3.7515,2.3602) -- (3.7664,2.5151) -- (3.7714,2.6714) -- (3.7664,2.8278) -- (3.7515,2.9827) -- (3.7266,3.1346) -- (3.6937,3.2758) -- (3.6499,3.4176) -- (3.5992,3.5466) -- (3.5401,3.6678) -- (3.4713,3.7826);
\draw[black,line width=.65pt] (3.4713,3.7826) -- (3.3944,3.8869) -- (3.3076,3.9815) -- (3.2175,4.0597) -- (3.1208,4.1249) -- (3.0179,4.1765) -- (2.9094,4.2138) -- (2.7958,4.2360) -- (2.6832,4.2429) -- (2.5670,4.2353) -- (2.4479,4.2132) -- (2.3265,4.1763) -- (2.2033,4.1248) -- (2.0791,4.0587) -- (1.9605,3.9825) -- (1.8364,3.8888) -- (1.7144,3.7826);
\draw[black,line width=.65pt] (1.7144,3.7826) -- (1.6039,3.6738) -- (1.4966,3.5561) -- (1.3868,3.4223) -- (1.2861,3.2866) -- (1.1848,3.1360) -- (1.0937,2.9864) -- (1.0083,2.8317) -- (0.9287,2.6714);
\draw[black,line width=.65pt] (0.9287,2.6714) -- (0.8303,2.4443) -- (0.7482,2.2152) -- (0.6875,1.9991) -- (0.6470,1.7925) -- (0.6295,1.6080) -- (0.6349,1.4428) -- (0.6463,1.3721) -- (0.6621,1.3130) -- (0.6854,1.2553) -- (0.7124,1.2094) -- (0.7451,1.1706) -- (0.7800,1.1421) -- (0.8242,1.1189) -- (0.8697,1.1055) -- (0.9207,1.1001) -- (0.9712,1.1025) -- (1.0327,1.1136) -- (1.0931,1.1317) -- (1.2284,1.1925) -- (1.3747,1.2826) -- (1.5399,1.4075) -- (1.7144,1.5603);
\draw[black,line width=.65pt] (1.7144,1.5603) -- (1.9304,1.7735) -- (2.1715,2.0373) -- (2.4222,2.3356) -- (2.6856,2.6714);
\draw[black,line width=.65pt] (2.6856,2.6714) -- (3.0681,3.1938) -- (3.4713,3.7826);
\draw[black,line width=.65pt] (3.4713,3.7826) -- (3.9523,4.5284) -- (4.4497,5.3429);
\draw[gray!50,thin] (0.0000,0.0000) rectangle (5.5000,5.3429);
\node[font=\small,fill=white,fill opacity=.85,text opacity=1,inner sep=1pt] at (1.0702,1.6126) {$S_{1}$};
\node[font=\small,fill=white,fill opacity=.85,text opacity=1,inner sep=1pt] at (2.6598,1.6953) {$S_{2}$};
\node[font=\small,fill=white,fill opacity=.85,text opacity=1,inner sep=1pt] at (1.7318,0.6913) {$S_{3}$};
\node[font=\small,fill=white,fill opacity=.85,text opacity=1,inner sep=1pt] at (1.7529,2.6736) {$S_{4}$};
\node[font=\small,fill=white,fill opacity=.85,text opacity=1,inner sep=1pt] at (3.2871,2.6693) {$S_{5}$};
\node[font=\small,fill=white,fill opacity=.85,text opacity=1,inner sep=1pt] at (4.5535,1.6331) {$S_{6}$};
\node[font=\small,fill=white,fill opacity=.85,text opacity=1,inner sep=1pt] at (2.6598,3.6476) {$S_{7}$};
\node[font=\small,fill=white,fill opacity=.85,text opacity=1,inner sep=1pt] at (1.0702,3.7303) {$S_{8}$};
\end{tikzpicture}
\par\smallskip (E8)\end{minipage}
\caption{Chebyshev curves of type $(3,4)$ (E7) and $(4,5)$ (E8).}\label{fig:fourth}
\end{figure}

\end{document}